# Assessing emission certification schemes for grid-connected hydrogen in Australia

Chengzhe Li, Lee V. White, Reza Fazeli, Anna Skobeleva, Michael Thomas, Shuang Wang, Fiona J. Beck*

## Abstract

Certification of low emissions products is crucial to ensure that they can attract a "green premium" in international markets but need to be designed so that they do not impose undue regulatory burden or drive-up production costs. In this work, we employ energy system modelling to evaluate how different policy choices affect the cost and certified emissions of hydrogen produced by grid-connected electrolysers in Australia. As a major exporter of energy and ores, Australia has ambitions to become a 'renewable energy superpower' by providing commodities with low embedded emissions. Our analysis shows that emissions are not always accurately accounted for under the currently proposed Australian methodologies. However, imposing geographic correlation and relatively relaxed temporal correlation requirements on the use of renewable energy certificates could be sufficient to ensure minimal emissions at low cost. Our findings provide evidence to assist in the design of robust and practical certification schemes.

## Introduction

Internationally recognized certification of embedded emissions is critical to enable trade in so-called "green" commodities like hydrogen and its derivatives where production processes and emissions can vary widely[1,2], allowing customers to differentiate between otherwise identical products[3,4]. In the absence of a global price on carbon, such differentiation will be needed for low-emission products to compete in global markets[5], and to facilitate the implementation of green industrial policies that incentivize decarbonization[6], for example, hydrogen production tax credits in the Inflation Reduction Act (IRA) of the United States scale inversely with certified emissions intensity[7].

In general, the certification of hydrogen produced via electrolysis is relatively straightforward when the power is sourced directly from dedicated, off-grid renewable energy (RE), such as wind or solar. Certification becomes more complex for grid-connected electrolyser systems, as the emissions from grid electricity provided by a dynamic mix of generators vary with time and location[8]. The challenge for certification schemes is to balance the requirement to accurately account for electricity emissions without imposing excessive regulatory burden, or driving up production costs[6].

To address this issue, the European Union (EU) has established a methodology in the Renewable Energy Directive II[9] (RED II) and associated Delegated Regulations[10,11] to certify electricity emissions associated with grid-connected hydrogen production under its CertifHy scheme[12]. Grid electricity is defined as "fully renewable", and hence zero emissions, if it meets three criteria: RE is supplied by renewable capacity commissioned no more than 36 months before the start of hydrogen production (additionality); RE should be generated and consumed by hydrogen production within a certain time interval (temporal correlation or simultaneity); and RE should be generated within the same interconnected bidding zone as the electrolyser (geographic correlation). If these conditions are not met, then the electricity emissions are calculated using the average emissions intensity of the electricity grid[10].

Among the RED II criteria, temporal correlation has proven to be an important and controversial policy lever[7,13]. While strict time-matching may lower the emissions of hydrogen by ensuring that the system can only use grid electricity when RE is being generated, it also impacts the profitability by limiting the operation of the electrolyser and requiring electricity or hydrogen storage to firm supply[14,15]. Recognizing the trade-off between economic competitiveness and minimizing emissions, jurisdictions have imposed less restrictive temporal correlation in the short term to enable grid-connected hydrogen production. For example, the EU's CertifHy has initially adopted temporal correlation over calendar months, and plans to gradually shift to hourly matching by 2030[11]. Likewise, the IRA in the US starts with the requirement of annual temporal correlation, gradually shifting to an hourly basis to assess qualification for the emissions intensity dependent tiers of the hydrogen production tax credit[16].

The Australian government is developing a Product Guarantee of Origin (PGO) scheme to certify emissions embedded in products such as hydrogen[17], to support its ambition to become "renewable energy superpower" by exporting low emission products and commodities[18,19]. A proposal for the methodology and emissions accounting framework that will underpin the PGO scheme has been released for consultation and the policy settings are under development[20–22]. Unlike the IRA or CertifHy, the Australian PGO scheme does not define emission thresholds for "low-emissions" or "renewable" hydrogen. Instead, the emission intensity for each batch of hydrogen is certified using a modular approach designed to meet the requirements of the International Partnership for Hydrogen and Fuel Cells in the Economy's methodology[23]. Accounting for emissions in each part of the supply chain allows interoperability with different international schemes with different definitions of "low-emissions hydrogen"[4].

The PGO scheme proposes two methods to estimate the emissions of electricity purchased from the grid, using methodologies defined under the National Greenhouse and Energy Reporting (NGER) scheme[24]. The location-based method applies government defined annual electricity emissions factors for each state in Australia, while the market-based method allows consumers of grid electricity to offset their emissions by purchasing and surrendering renewable energy certificates (RECs) which are generated for each megawatt hour of renewable electricity sold to the grid. Each REC will report the commissioning date and location of the power station, and a hourly time stamp recording when the electricity is generated[17]. While the information included with each REC grants the PGO scheme flexibility to develop rules to ensure additionality as well as geographic and temporal correlation, optimum settings have not yet been determined.

The Australian National Electricity Market (NEM) is both a wholesale market and a physical power system, providing electricity to 88% of Australia's population living along the eastern seaboard[25]. The NEM is one of the longest transmission networks in the world, connecting five state grids which each have unique mix of electricity generators[25]. South Australia has rapidly increased the penetration of renewables on its grid over the last 20 years, reaching 71% in 2023[26], and is on track to be powered by 100% renewable energy by 2027[27,28]. Queensland has exceptional solar resources, but started the transition much later, reaching 28% in 2023 and aiming for 80% by 2035[29]. The different levels of ambition present unique challenges for designing effective national policy settings for the PGO scheme. While previous analyses have shown the importance of geographic and temporal correlation requirements on costs and emissions in European and US grids[7,13–15,30], there has been no analysis of the impact of these policy choices on the competitiveness of hydrogen produced across states in Australia.

In this context, our study aims to assess the implications of different policy settings in the proposed PGO certification scheme on the cost and emissions intensity of hydrogen produced by grid-connected hydrogen production systems in Australia. We establish a linear optimization model to quantify the difference between the actual emissions and those assigned by the proposed PGO accounting methods for different cost-optimal system configurations. We explore how and why proposed emissions accounting frameworks under the PGO will lead to inaccuracies and identify effective and efficient policy settings that can keep both emissions and costs low. Lastly, we provide policy suggestions for the future development of the Australian PGO scheme.

# Hybrid grid-connected hydrogen production

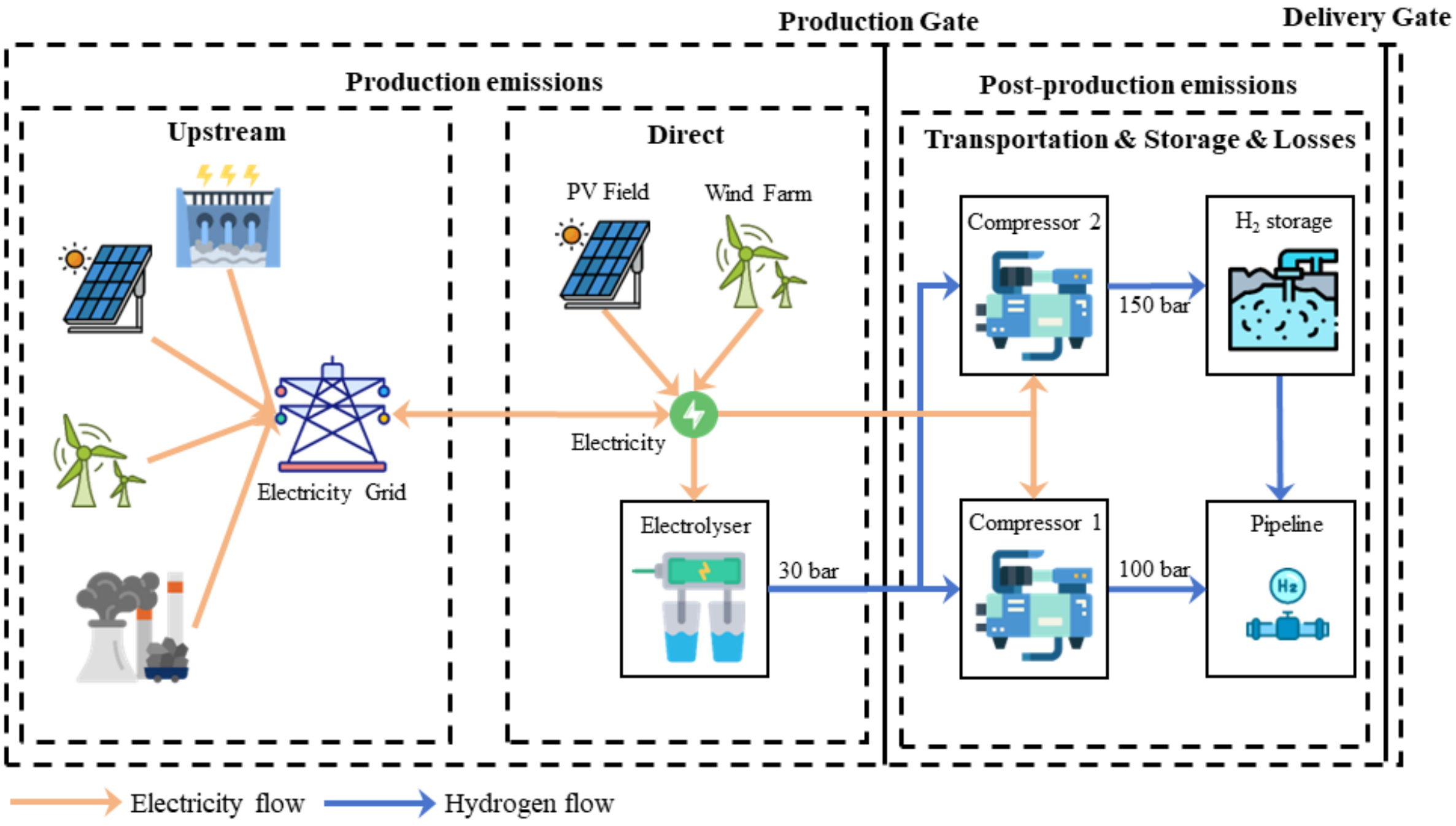

**Figure 1.** System boundaries for emissions accounting of hybrid grid-connected hydrogen production under the Product Guarantee of Origin scheme. The arrows indicate energy commodity flows and their directions: orange arrows for electricity, and blue arrows for hydrogen.

The emission accounting system boundary of the Australian PGO scheme split into two main parts[20]: 1) the production boundary encompassing upstream and direct emissions; and 2) the post-production boundary including emissions from transport and storage of hydrogen. The PGO scheme requires that hydrogen producers define a batch period over which they are certifying their hydrogen production, ranging from hourly to yearly. During this batch period, the producer collects emissions and RECs information and reports it to the PGO scheme to label the certified emission intensity for the batch of hydrogen.

In this study, we evaluate the LCOH and associated actual and certified emissions for an optimized grid-connected hydrogen production system, powered by a mix of renewable and grid electricity as shown in **Figure 1**. Production emissions primarily stem from upstream scope 2 greenhouse gas (GHG) emissions attributable to electricity imports from the grid. Emissions from the import of potable water are included in the PGO scheme[20] but are not included in our model for simplicity, as they would be identical whether hydrogen production is powered by off-grid RE or on-grid electricity. In the post-production boundary, we include scope 2 emissions associated with electricity used for compression, transport and storage. We neglect fugitive emissions of hydrogen which has a GWP of 11[31] which could occur at various points during the supply chain.

We establish an energy system model of the grid-connected hydrogen production system, based on previously developed models of off-grid systems[32] consisting of dedicated wind and solar, an electrolyser, two compressors, hydrogen storage facilities and pipeline (details in **Method** section). We assume that the system can only access RECs from RE capacity installed at the same time as the electrolyser, which can be directly connected, or accessed through the grid. This assumption allows us to accurately track emissions associated with RECs and reflects the rapid deployment of new capacity to meet the growing demand for RECs in Australia[33].

The system is operated with hourly time resolution and delivers a constant supply of hydrogen to industrial end-users throughout the year[34]. The system can use power directly from the dedicated RE capacity or can choose to buy or sell electricity from the Australian NEM, modelled using historical data on spot prices and emission factors for the five interconnected states[35]. Purchasing electricity from the grid incurs a "transmission use of system fee" representing the costs required to maintain the electricity grid set by grid operators[36,37]. A linear optimization model is used to determine the capacity of the different components as well as the interaction with the grid required to ensure the lowest hydrogen production cost, calculated as the levelized cost (LCOH). The capacity of the system is constrained so as not to perturb the spot price or emission intensity of the grid.

We leverage our model to assign time and state dependent emission factors to net grid electricity consumption in each time step. Average emissions factors (AEFs) are computed by averaging emissions generated by all electricity generated across each of the five states, and marginal emissions factors (MEFs) represent the emissions of the marginal generator, which is the generator that would be dispatched to meet additional electricity demand. MEFs best represent the actual electricity emissions caused by the additional demand on the grid from our system[38–40] and quantify the grid emissions avoided by displacing electricity from the marginal generator with RE sold into the grid. Comparing the emissions calculated using the PGO methods described previously with the MEF-based method allows us to quantify how well the proposed certification methodology captures the actual emissions.

Policy settings can be incorporated into the model by introducing additional constraints. Geographic correlation can be imposed by requiring the dedicated RE capacity to be in the same state as the electrolyser. Temporal correlation requires that the total electricity sold to the grid must be larger than or equal to the total purchased electricity within the time interval[14].

# Evaluating certified emission intensity of hydrogen

[34] Note that the constant hourly supply is to ensure the hydrogen supply reliability required by refineries, ammonia plants, liquefaction facilities, and other industrial operations.

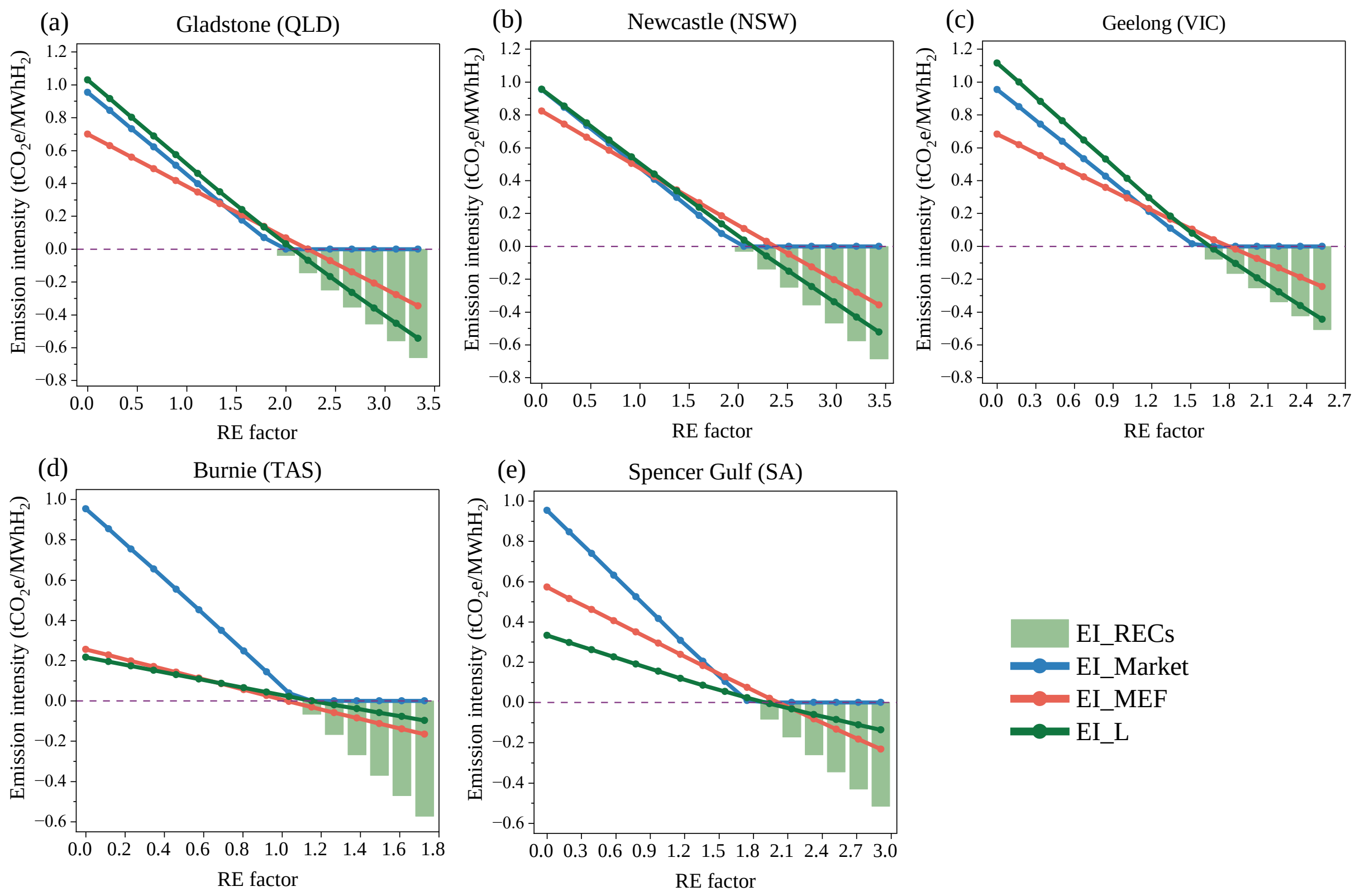

**Figure 2 (a)-(e)** Greenhouse gas **(**GHG) emission intensities of hydrogen under various RE factors, defined as the ratio between RE generators and electrolyser capacity in five states of the NEM, calculated using three different methods for electricity emissions. Actual emissions are tracked using time and state dependent marginal emission factors to net grid electricity consumption (EI_MEF, orange curve). Certified emissions are calculated using the market-based (EI_Market, blue curve) and location-based (EI_L, green curve) adopted in the PGO scheme. Negative emissions intensity values indicate that more renewable electricity is sold to the grid than is bought per unit hydrogen produced, resulting in a reduction of grid emissions (for location and MEF methods). When the system generates more RECs that results EI_Market = 0[24] and grid electricity emissions can be reduced by excess RECs generated for each MWh hydrogen produced (EI_RECs, green bars). The conversion rate between energy content units of hydrogen throughout the manuscript is based on its higher heating value (39.4 kWh/kg $H_2$[41]).

We calculate the emission intensity of hydrogen produced by hybrid systems in the five states for systems that maintain a constant electrolyser size (equal to that of an optimized off-grid system, around 10-50MW) but have varying levels of dedicated RE capacity (**Figure 2)**. The dedicated RE capacity is co-located ensuring geographic correlation and the systems can interact with the grid flexibly without temporal correlation. The RE factor is defined as the ratio between RE generators (sum of wind and solar) and electrolyser capacity: an RE factor of 0 is a fully grid-reliant system, while an RE factor of 2 indicates the system has RE capacity two times larger than the electrolyser. For each location, the largest RE factor is set to be 1.5 times greater than the RE factor for the optimised off-grid system. By adjusting the RE factor from 0 to its maximum value, we explore the variation in the emission intensities of hydrogen using the three different accounting methods.

In general, comparison between the emission intensities calculated using the AEF-based and PGO location-based methods, given in **Figure 15**, show relatively good agreement for all systems investigated, as they both use state-based AEFs, albeit at different temporal resolutions (1 hour for EI_AEF, 1 year for EI_L)[42]. However, the location-based method does not reproduce the actual emissions intensity well, calculated with MEF-based method. It overestimates actual emissions in QLD, NSW, and VIC by up to 63% and underestimates actual emissions for SA by 42% (RE Factor = 0). This is because MEFs for states dominated by coal (QLD, VIC, NSW) tend to be lower than the AEFs used in the location-based method, as over 49% of the time the marginal generator will be a lower emissions technology like gas or renewables (**Figure 13**). Conversely, SA tends to have MEFs higher than AEFs as coal is more likely to be the marginal generator when importing electricity.

The PGO market-based method overestimates the emissions for fully grid reliant systems (RE Factor = 0) and overestimates the electricity emissions avoided by self-sufficient systems (maximum RE Factor) by generating excess RECs. The largest discrepancies are observed for TAS and SA, with emissions intensities overestimated by up to 274%, and EI_RECs overestimated by up to 245% in TAS. Deviations reduce as less grid-electricity is bought and sold.

We attribute these discrepancies to the adoption of nation-wide factors in the PGO market-based method (**Eq. 18a**). The proportion of RE on the state grids can be much higher (93% and 71% in 2023 for TAS and SA respectively, see **Table 7**) than the national-wide applicable renewable power percentage (ARPP) of 18.72%. Likewise, the emissions intensity used in the calculation (RMF) is calculated by adjusting the average emission factor of the grid at the national level (including the NEM, South-West Interconnected System, and the Darwin, Katherine Interconnected System[43]) by excluding all RE that has generated RECs (**Eq. (27)**). This overlooks the energy mix of regional grids, resulting in larger discrepancies for states with higher RE penetration, overestimating emissions when buying grid electricity.

# Impact of temporal correlation

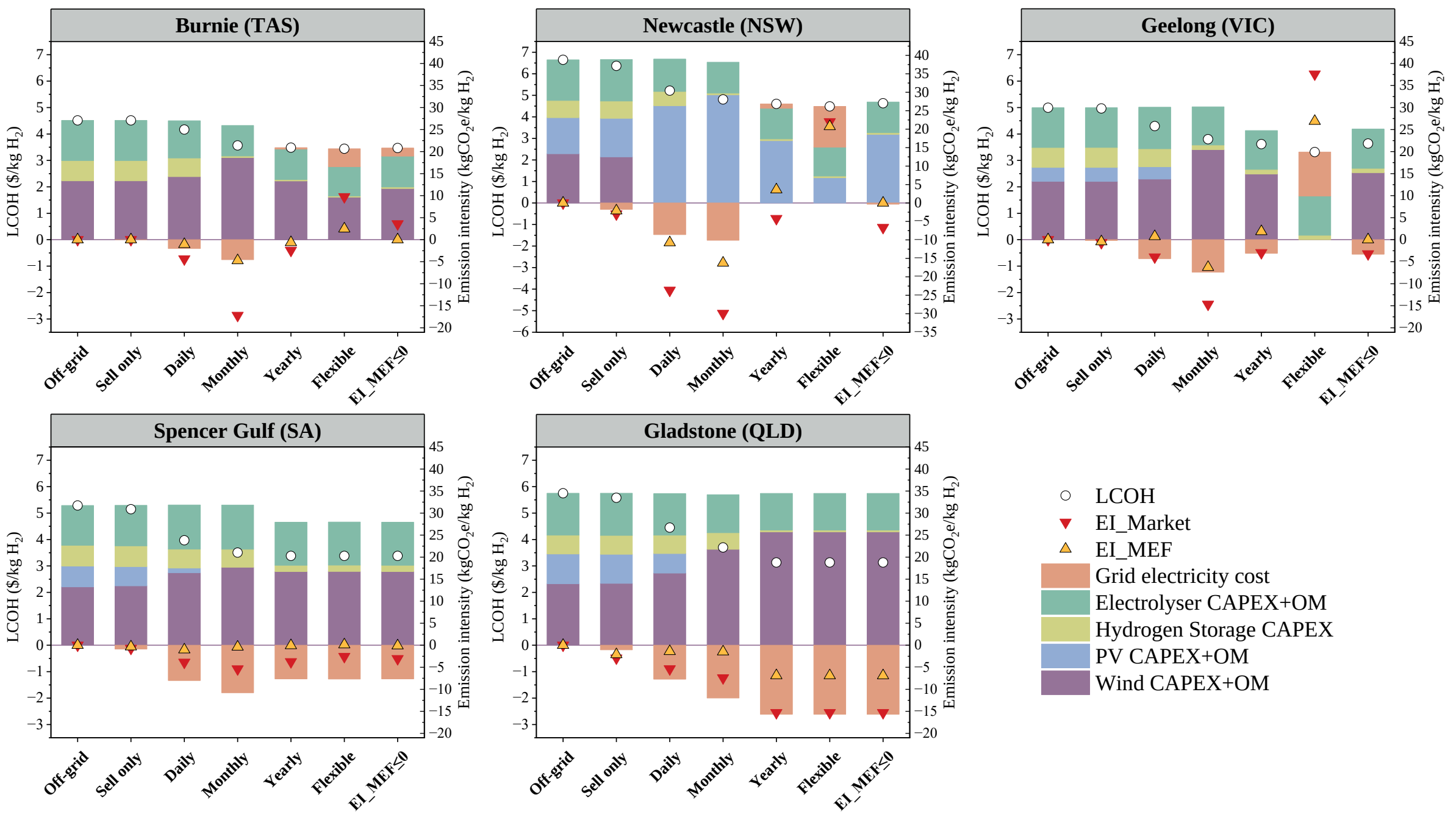

**Figure 3** Total LCOH (white dots), average yearly actual (EI_MEF, yellow triangles) and certified emission intensities of hydrogen (EI_Market, red triangles) calculated using the PGO market-based method, optimized system configurations in five states across different under different constraints described in **Table 1**. The emission intensity is given in kgCO$_2$e/kgH$_2$ to allow comparison with existing certification schemes. Bars show LCOH deconvoluted to components: grid electricity cost (sum of the total expenditure and profit from buying and selling electricity to and from the grid, which can be negative), capital expenditure (CAPEX) of hydrogen storage, and CAPEX and operation and maintenance costs (OM) of electrolyser, PV and wind generators. Negative emission intensity values for market-based method (EI_Market) quantify grid electricity emissions offset by excess RECs generated for each kg hydrogen produced, while they represent the actual emissions displaced for MEF-based method (EI_MEF).

**Table 1** Scenario design

| Scenario | Grid connection | Temporal correlation | Emission intensity constraint |
|---|---|---|---|
| Off-grid | NO | None | None |
| Sell only | Yes – but only sell electricity | Hourly | None |
| Daily | Yes – buy and sell | Daily | None |
| Monthly | Yes – buy and sell | Monthly | None |
| Yearly | Yes – buy and sell | Yearly | None |
| Flexible | Yes – buy and sell | None | None |
| EI_MEF≤0 | Yes – buy and sell | None | $e_{MEF} \leq 0$ |

In the following, we use our model to optimize the configuration and operation of hybrid grid-connected hydrogen systems to minimize LCOH under constraints, as described in **Table 1**. We compare the scenarios include the reference off-grid scenario where no grid interaction is permitted; four scenarios requiring temporal correlation over different intervals (hourly - which is equivalent to only allowing the selling of electricity to the grid labelled "Sell only", daily, monthly, yearly), while ensuring geographic correlation by co-locating the dedicated RE capacity, and a "flexible" scenario with no requirements for temporal correlation. To assess the effectiveness of temporal correlation policies for lowering emissions, we also include a scenario where the system can interact with the grid without temporal correlation, and instead constrain the "actual" emissions intensity of hydrogen to be less than zero, i.e. EI_MEF≤0. This final scenario is not meant to model a realistic policy setting, but to understand what a system optimised for zero-emissions at lowest cost looks like. We focus on the PGO market-based method, that allows the trading of RECs and explore the impact of different temporal correlation policy settings on reducing actual emissions, and those certified by the PGO scheme.

LCOH is highest for off-grid production systems and reduces with grid connection for all states, with the lowest LCOH occurring for fully flexible scenarios with no additional constraints **(Figure 3)**. There are two distinct trends for optimized system configurations and emission intensity of hydrogen: one that holds for TAS, NSW and VIC, and one that holds for SA and QLD.

Due to the interplay between spot prices and the quality of the RE resources, the cheapest way to produce electrolytic hydrogen in TAS, NSW and VIC is to minimize capital investment and source large amounts of power from the grid. For TAS and VIC, relatively low grid electricity prices make dedicated RE uncompetitive. While NSW has the highest yearly average electricity price in the NEM, the local RE resources in Newcastle are less abundant (quantified by the RE capacity factor, **Table 6**), leading to higher RE costs. When unconstrained by temporal correlation in the flexible scenario, the system builds less dedicated RE capacity in TAS and NSW compared to the off-grid scenario, and none in VIC, and invests in smaller hydrogen storage and electrolyser capacity. Critically, relying on grid power in the lowest cost flexible scenarios can result in very high hydrogen emission intensities, with EI_MEF =~25 $kgCO_2e/kgH_2$ for VIC and NSW, comparable to hydrogen made from coal gasification[44]. Once again, we see that the PGO market-based method overestimates the certified emissions intensity, particularly for TAS and VIC.

Once yearly temporal correlation is imposed in TAS, NSW and VIC, the system is forced to build RE capacity, the LCOH increases, and the actual emission intensity falls to below zero in TAS and slightly above zero in NSW and VIC. As the temporal correlation requirements become more restrictive, shifting from yearly to hourly time matching, the system increases the capacity of the electrolyser and storage, and the configuration progressively resembles an off-grid setup. Consequently, LCOH increases by 30%-38%. Once again, the PGO market-based emission calculation overestimates potential emissions offsets by RECs, becoming more significant when the system sells larger amounts of electricity to the grid.

The trend is different for SA and QLD. All systems have actual emissions intensities close to or below zero, and the flexible and yearly temporal correlation scenarios all result in nearly identical system configurations, with similar and lowest cost LCOH. The cost reduction in these two scenarios compared to the off-grid system is mainly attributable to the systems' ability to profit from selling surplus electricity to the grid, driven by favourable local RE resources and higher regional electricity prices (**Table 6, 7**). The system significantly oversizes the cheapest RE capacity and sells large amounts of electricity to subsidize the cost of hydrogen production. As the temporal correlation requirements become more restrictive, shifting from yearly to hourly time matching, the restrictions on importing electricity from the grid results in the establishment of configurations similar to the off-grid system with PV and larger storage components. The systems are prevented from subsidizing the cost of hydrogen by selling electricity, increasing the LCOH by 78 % in QLD, and 52% in SA.

One way to assess the efficacy of temporal correlation as a policy is to compare with systems optimized to minimize cost while constraining $EI_MEF \leq 0$. Figure 3 shows that such systems are most similar to those optimized under yearly temporal correlation, suggesting it is an efficient policy measure to keep emissions low, at low cost. However, it is important to note that while yearly temporal correlation is sufficient to ensure that the certified emissions intensities under the PGO market-based method are all below zero, actual emissions intensities are still 1.9 kgCO2e/kgH2 for the VIC case and 3.6 kgCO2e/kgH2 for the NSW case. These emissions intensities are similar to those of natural gas based blue hydrogen systems[1].

# Impact of geographic correlation

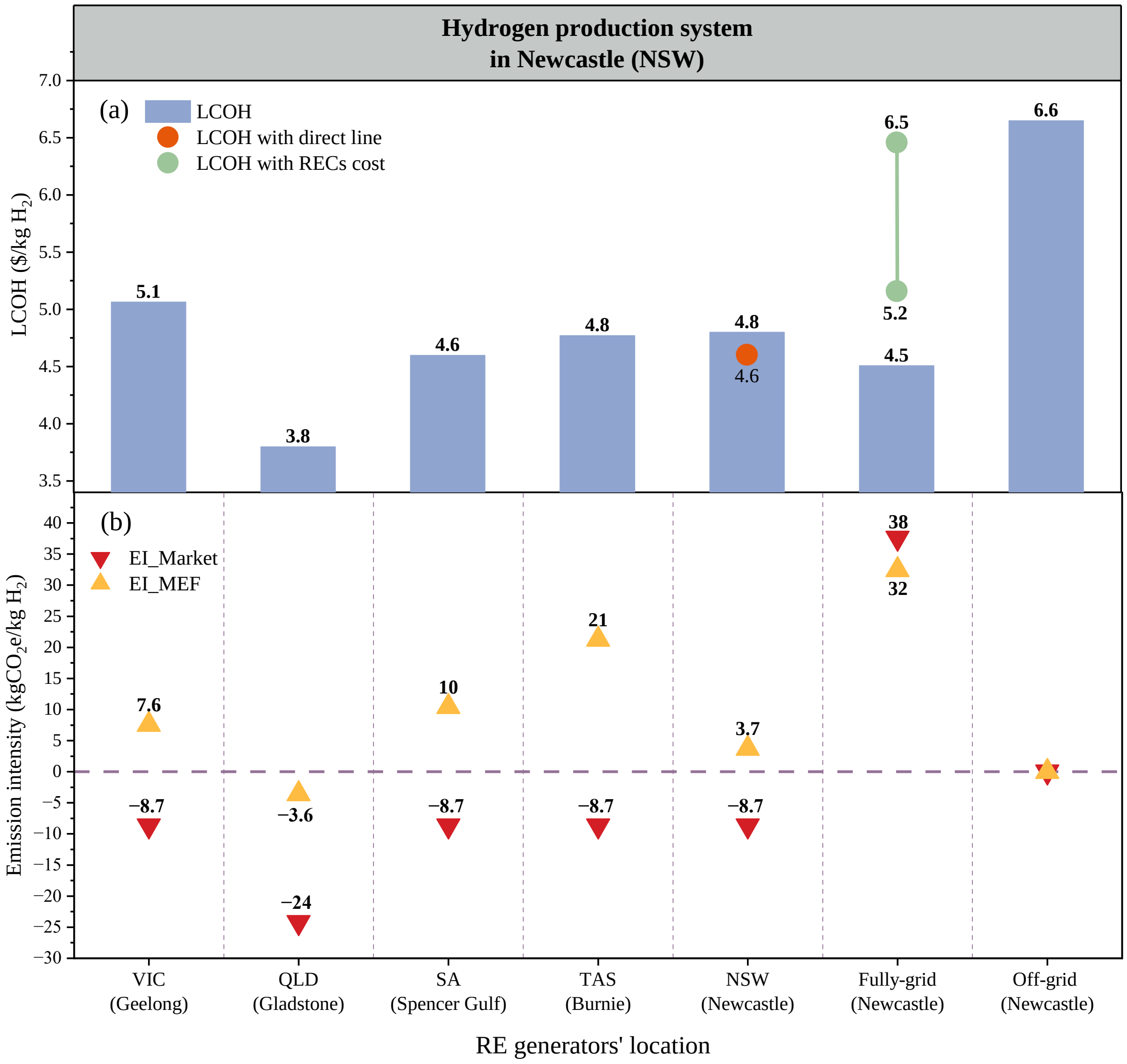

**Figure 4 (a)** LCOH and **(b)** emission intensity of hydrogen produced by a grid-connected electrolyser located in Newcastle NSW, with RECs sourced from grid-connected RE capacity located different states, as indicated along the x-axis. Systems are optimised under the constraint of yearly temporal correlation and are compared to off-grid and fully-grid connected systems in NSW. The green line indicates the estimated cost of purchasing RECs from a third party to offset emissions for the fully on-grid case, assuming REC costs of 20-60AUD[33]. Note that the LCOH for a system with RE located in NSW connected to the electrolyser through the grid is slightly higher than that calculated previously (see **Figure 3**) for a system with direct connection between co-located RE and electrolyser in NSW (red dot) due to increased transmission use of system fees. The emission intensity of hydrogen is calculated by MEF-based method (yellow triangles) and for the PGO market-based (red triangles).

So far, all hybrid systems have been constrained by geographic correlation, with a direct connection between co-located RE and hydrogen electrolysers. Next, we relax that constraint and compare scenarios in which the electrolyser and the dedicated RE generators are located in different states and interact only through grid-connection, as would occur under certain pre-purchase agreement[39]. In these scenarios, the RE is sold to its local grid, generating RECs and offsetting actual emissions in that state, while the electrolyser buys electricity from its local grid, generating actual emissions and using the RECs generated by the RE to offset emissions (see Modelling geographic correlation in the **Method** section). We focus on an electrolyser located in Newcastle (NSW), which is the site of the announced Hunter Valley Hydrogen Hub, slated to include a 50 MW electrolyser powered by grid-connected RE through purchase and retirement of RECs[45]. We compare systems under the constraint of yearly temporal correlation with dedicated RE located in different states, with a fully grid-powered system and an off-grid system (**Figure 4**). The LCOH can be reduced by as much as 21% when the RE is located in QLD due to the superior RE resources and high electricity costs, as the net cost of grid electricity is determined by buying and selling electricity to and from two different grids and depends on the spot prices in both states.

Results show that yearly temporal correlation is sufficient to ensure that all certified emissions (EI_Market, red triangles) are below zero even without geographic correlation, as each unit of electricity used by the electrolyser is compensated for by a REC generated within the year. However, the actual emissions (EI_MEF, yellow triangles) vary widely with the location of the RE as the emissions displaced by selling RE will depend on the makeup of its local grid and are not necessarily equal to emissions generated using electricity in NSW. In the worst case, the emission intensity can reach 21 $kgCO_2e/kgH_2$, similar to that of hydrogen produced via coal gasification[44], as high emissions intensity electricity consumed in NSW is being offset by displacing relatively low emissions intensity electricity in TAS.

In all scenarios so far, we assume a strict additionality requirement: that dedicated RE is deployed concurrently to the electrolyser. The PGO scheme does not require such strict additionality as RECs can be generated by eligible RE systems built before the electrolyser, as long as they began generating renewable electricity after 1997, and do not have to be surrendered as part of Renewable Energy Target[46–48] or any other scheme. The fully grid-powered electrolyser scenario in **Figure 4** approximates this case, buying electricity from the NSW grid at the spot price and generating the corresponding emissions according to the state-based MEF. Under this scenario, the additional cost for RECs to offset the certified emissions under the market-based method would increase the LCOH by 16-44% assuming RECs costs of AUD 20-60[33], while estimating the actual emissions would require the time and state dependent MEF for each REC.

[45] The LCOH of the NSW system connected through the grid is slightly higher than when it is directly connected with RE, as plotted in **Figure 3**, since electricity purchases from the grid and incurs a “transmission use of system fee” (see **Table 3**)

# Conclusion

Australia has the potential and the ambition to be a producer and exporter of large volumes of low-emissions hydrogen. This paper evaluates the efficacy of proposed emissions accounting methodologies and certification scheme rules for hydrogen production in Australia's PGO certification scheme, which is currently under development.

We employ an energy system model to quantify emissions and costs of hybrid grid-connected hydrogen production systems located across Australian NEM. By imposing different constraints, we assess the accuracy of the PGO certification scheme methodologies and the efficacy of different policy settings in minimizing emissions while maintaining low costs.

We find that PGO emissions accounting methods can lead to significant discrepancies compared to the actual emissions. The differences vary across states within the NEM due to differences in state grid generator mix and renewable resource availability. The PGO location-based method overestimates emissions when AEFs are larger than MEFs (e.g. 63% for fully grid reliant in VIC), and underestimates when they are smaller (e.g. 42% for fully grid reliant in SA). This occurs because actual emissions are determined by the marginal generator, not the average grid emissions intensity. The high frequency of hydro as a marginal generator leads to overestimation in heavily coal-reliant states with high AEFs (QLD, NSW, VIC), while coal results in underestimation in high renewable penetration states with low AEFs (SA).

The PGO market-based method tends to overestimate emissions in all states when large amounts of electricity is bought (e.g. 274% for fully grid-reliant systems in TAS). The discrepancies are larger for states with more RE penetration as the market-based method uses factors relating to RE penetration and residual emissions (i.e. ARPP and RMF) that are currently estimated using the RECs information aggregated at the national level[43]. It is likely that the use of more granular, state-level factors would improve accuracy.

We note that the overestimation is consistent with the principle of conservativeness recommended for embedded emissions accounting[6], and hence is not necessarily a design flaw, underestimation of embedded emissions contravenes the conservativeness principle and undermines the scheme's emissions-reduction goals.

We further investigate the effect of PGO policy settings on the cost as well as the real and certified emissions when trading renewable energy certificates (i.e. PGO market-based method). By ensuring that the systems only use RE from dedicated, co-located wind and solar capacity, we are implicitly applying additionality and geographic correlation criteria. We then examine different temporal correlation settings. Despite differences in cost optimal system configurations across the five states, our main finding is robust and in line with previous work[7,13]: allowing less onerous temporal correlation requirements can reduce LCOH, while still limiting emission intensities. Specifically, optimizing systems operating under yearly temporal correlation results in similar system configurations as those optimized to keep certified emissions below zero at lowest cost (i.e. constraint of EI_MEF≤0). This is an important outcome that suggests that less stringent temporal correlation is an effective and efficient policy measure to ensure near-zero emissions systems at low cost and low regulatory overhead. However, in certain circumstances, discrepancies between the actual and certified emissions calculated under the PGO market-based method could lead to systems that still generate hydrogen with actual emissions under yearly correlation, even though they are certified as zero emissions, particularly when the geographic correlation is not met (e.g. 21 $kgCO_2e/kgH_2$ in **Figure 4**). Ensuring actual emissions are below zero for all states will require geographic correlation coupled with temporal correlation intervals somewhere between monthly and yearly and will change with location and over time, as it is determined by the interplay between regional electricity prices, local renewable resources and RE penetration in the local state grid.

Our findings provide timely evidence to aid the design of the Australian PGO scheme as well as evidence to support the development of robust and practical certification schemes internationally.

# Method

## Model description

To quantify the cost and emission intensity of hydrogen production by the hybrid grid-connected hydrogen production system, we establish a linear optimization model to determine the optimal capacity of the system components consisting of dedicated wind and solar, an electrolyser, two compressors, hydrogen storage facilities and pipeline, and optimal interaction with the grid to ensure the lowest hydrogen production cost. The model is operated with hourly time resolution over a one-year period and is constrained to deliver a constant supply of hydrogen to industrial end-users throughout the year. The modelled system can use the dedicated RE capacity to supply power for electrolysis and subsequent compression processes, (directly or through the grid) or can choose to buy or sell electricity from the Australian NEM. The system can exchange unlimited power with the grid however RE generation and electrolyser capacities are kept within the range of 10MW to 50MW so that the effect on the price and emission intensity of grid electricity can be neglected. These values were chosen to reflect the capacity of proposed hydrogen production projects in Australia[49], and the average instantaneous electricity demand of the smallest regional grid in the NEM, Tasmania, which is approximately 1200MW[50]. Local wind and solar electricity generation are calculated using NREL's PySAM engine[51] based on historical local weather data from NASA's MERRA-2 dataset[52]. The electrolyser converts electricity into hydrogen with an efficiency of 70%[53]. We assume that the hydrogen produced by the electrolyser is initially at a pressure of 30 bar and compressor 1 pressurizes hydrogen to 100 bar for direct pipeline transport to the load, while compressor 2 increases the pressure to 150 bar for storage in underground lined rock caverns (LRC) if needed. The $H_2$ storage acts as a buffer to the load, allowing it to be stockpiled when local RE generation is abundant and grid electricity is cheap, and to be supplied when local RE generation is scarce, and electricity is expensive. The system could choose to store electricity in batteries; however, due to the high capital expenditure of batteries compared to hydrogen storage, the model optimization process favours investing in hydrogen storage. The system can also choose to import and export electricity from the state grid (bidding zone) in which it is located, defined by time and location dependent historical electricity prices and emissions factors, described in detail in section **Modelling the Australian National Electricity Market** below.

## Hybrid grid-connected hydrogen production optimization model

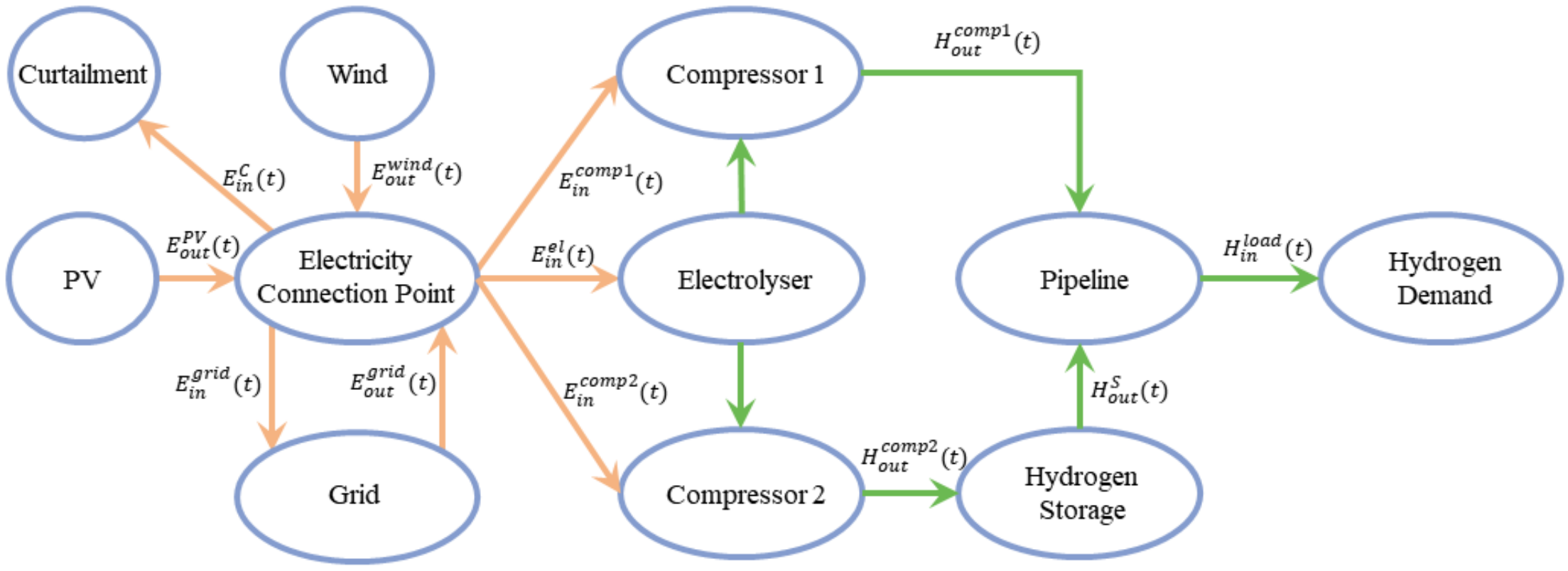

**Figure 5** Overview of the hybrid grid-connected hydrogen production optimization model. Ovals represent nodes, which are system components responsible for transformation of energy commodities. Lines represent edges which indicate the flow of an energy commodity and arrows represent the flow directions. Orange lines indicate electricity flow (kW), while green lines represent hydrogen flow (kg/hour). The electricity connection point serves as the node ensuring electricity balance, where electricity consumption aligns with electricity generation. The grid node represents the Australian NEM, and the hydrogen production system can participate in the Australian NEM via the transmission line, either as a generator selling electricity or as a consumer purchasing electricity.

The linear optimization model is shown schematically in **Figure 5**, where ovals correspond to nodes representing system components for converting between energy commodities, lines representing edges which indicate the energy commodity flows, with arrows indicating the direction of these flows. The variables $E_y^x$ and $H_y^x$ indicate flows of electricity in kilowatts and hydrogen in kilograms per hour between different component, respectively, where the superscript $x$ specifies the component, and $y$ indicates the flow direction. The production system simulates operations by ensuring the equilibrium of electricity and hydrogen flows, as described in the following.

### *Electricity flow balance*

The production system is powered by electricity which can be sourced from local RE generation or imported from the grid. At each time $t$, local wind generation ($E_{out}^{wind}(t)$) and solar generation ($E_{out}^{PV}(t)$) are produced based on the local weather profile. Renewable energy (RE) can either be used to (i) power the electrolyser for hydrogen production ($E_{in}^{el}(t)$) and the compressor for hydrogen transmission into either the pipeline ($E_{in}^{comp1}(t)$) or storage facilities ($E_{in}^{comp2}(t)$), (ii) exported to the grid ($E_{in}^{grid}(t)$), or (iii) curtailed ($E_{in}^{C}(t)$), where the variables refer to power flows shown in **Figure 5**. If the local RE is insufficient to support the system's operation, the system can choose to import electricity from the grid ($E_{out}^{grid}(t)$). The decision of electricity dispatching is made by solving the flow balance equation at the Electricity Connection Point given by

$$E_{in}^{el}(t) + E_{in}^{comp1}(t) + E_{in}^{comp2}(t) + E_{in}^{grid}(t) + E_{in}^{C}(t) = E_{out}^{wind}(t) + E_{out}^{PV}(t) + E_{out}^{grid}(t) \quad \forall t \quad (1)$$

### *Hydrogen flow balance*

The system converts electricity into hydrogen through electrolysis and then either pumps the hydrogen into the pipeline for direct transportation to the hydrogen load or stores it in the hydrogen storage equipment, which serves as a reserve for the load, storing hydrogen to supply when production by the electrolyser is insufficient to meet the demand. At each time $t$, hydrogen is produced by the electrolyser with efficiency, $\eta = 70\%$[53], as following:

$$H_{out}^{el}(t) = \frac{E_{in}^{el}(t) \times \eta}{HHV} \qquad \forall t \quad (2)$$

where $HHV$ =39.4 kWh/kg$H_2$[41] is the higher heating value of hydrogen and $H_{out}^{el}(t)$ represents the hydrogen produced by electrolyser. This hydrogen is subsequently allocated by the system to either compressor 1 ($H_{out}^{comp1}(t)$), which pumps it into the pipeline for direct transportation to end-users at a pressure of 100 bar, or to compressor 2 ($H_{out}^{comp2}(t)$), where it is pressurized to 150 bar and then stored in the hydrogen storage:

$$H_{out}^{el}(t) = H_{out}^{comp1}(t) + H_{out}^{comp2}(t) \qquad \forall t \quad (3)$$

The hydrogen load in each time $t$ ($H_{in}^{Load}(t) = 180kg$) is met either by hydrogen flow directly from the electrolyser via pipeline and/or hydrogen from storage ($H_{out}^{S}(t)$) according to:

$$H_{in}^{Load}(t) = H_{pipe}^{comp1}(t) + H_{out}^{S}(t) \qquad \forall t \quad (4)$$

### *Wind and solar generation*

Local wind and PV generation at each time $t$ depends on the optimized capacity of the generators. To enable the optimization, we define reference capacities for the wind farm ($C_{Ref}^{wind} = 320$ MW)[54] and PV field ($C_{Ref}^{PV} = 1$ MW), and calculate the hourly reference generation ($E_{Ref_{out}}^{wind}$ and $E_{Ref_{out}}^{PV}(t)$) using NREL's PySAM engine[51] and historical local weather data from MERRA-2 dataset[52]. The actual generation is then calculated through linearly scaling the reference output based on the ratio of actual capacity to reference capacity:

$$E_{out}^{wind}(t) = \frac{C^{wind}}{C_{Ref}^{wind}} \times E_{Ref_{out}}^{wind}(t) \qquad \forall t \tag{5}$$

$$E_{out}^{PV}(t) = \frac{C^{PV}}{C_{Ref}^{PV}} \times E_{Ref_{out}}^{PV}(t) \qquad \forall t \tag{6}$$

Where $E_{out}^{wind}(t)$ and $E_{out}^{PV}(t)$ represent the actual hourly generation of the wind farm and the PV field, respectively. $C^{wind}$ and $C^{PV}$ indicate the capacity for the wind farm and PV field, respectively.

### *Grid node*

We introduce the grid node to simulate the interaction between production system and the grid. The grid is defined by historical electricity prices (spot prices) and emissions factors (average and marginal) data sourced from the NEM. The modelling approach is introduced in the following section **Modelling the Australian National Electricity Market**. Detailed information on the NEM data is provided in following section **NEM Data**.

### *Compressor 1 and compressor 2*

The electricity consumption of compressor 1, which pressurizes hydrogen produced by the electrolyser for pipeline injection, and compressor 2, which pressurizes hydrogen for storage is given by

$$E_{in}^{comp1}(t) = H_{out}^{comp1}(t) \times \mu_{out}^{comp1} \qquad \forall t, \tag{7}$$

$$E_{in}^{comp2}(t) = H_{out}^{comp2}(t) \times \mu_{out}^{comp2} \qquad \forall t, \tag{8}$$

where $\mu_{out}^{comp1}$ and $\mu_{out}^{comp2}$ are the electricity consumption for each unit of hydrogen compressed by compressor 1 and compressor 2, respectively.

### *Hydrogen storage*

The hydrogen storage level at each time $t$ ($H^{S_level}(t)$) is updated according to:

---

[54] Note that a wind farm requires a layout of multiple turbines to account for the wake effect between them. Here, we assume that 320 MW capacity represents a sufficiently large plant to consider a turbine layout.

$$H^{S_{level}}(t) = H^{S_{level_0}} + \sum_0^t \left(-H_{out}^S(t) + H_{out}^{comp2}(t)\right) \qquad \forall t, \tag{9}$$

where $H_{out}^S(t)$ represents the hydrogen sent out by the hydrogen storage to meet the load and $H^{S_level_0}$ is the initial level of hydrogen storage as we assume the system has been in trial operation for a period. The hydrogen storage level in the last time point should be equal to the initial level to ensure that all the hydrogen sent out is produced within the modelled year.

***Optimization***

The objective of the optimization model is to minimize the LCOH, defined as[55]:

$$LCOH = \frac{CAPEX \times CRF + O\&M + C^e}{M_{H_2}}, \tag{10}$$

where $M_{H_2}$ is the total hydrogen produced annually, $CAPEX$ is the capital expenditure including installation and equipment cost[53], O&M is the operation and maintenance costs, $CRF$ is the capital recovery factor, and $C^e$ is the electricity cost over the year. The CAPEX and O&M are given by the sum of the costs associated with individual system components:

$$CAPEX = \sum_{k \in K} C^k I^k, \tag{11}$$

$$O\&M = \sum_{k \in K} C^k \left(FOM^k + VOM^k\right), \tag{12}$$

where $C^k$ represents the installed capacity of a component, $I^k$ is investment needed per unit installed capacity, $FOM^k$ and $VOM^k$ are the fixed and variable operation and maintenance cost of component $k$, where $K$ is the component set which includes wind, PV, electrolyser and $H_2$ storage facilities. The $CRF$ can be calculated according to the plant lifetime $n$ and the interest rate $i$ as:

$$CRF = \frac{i(1+i)^n}{(1+i)^n - 1}, \tag{13}$$

The electricity cost over one year includes costs associated with buying electricity from the grid and negative "costs" from selling electricity to the grid in each time $t$ and is given by,

$$C^e = \sum_{t=0}^{8759} (E_{out}^{grid}(t) \times (P_j(t) + TS) - E_{in}^{grid}(t) \times P_j(t)), \tag{14}$$

where $P_j(t)$ is the electricity spot price in state $j$ at each time $t$. In addition, the system will need to pay an additional transmission use of system fee $TS$ when importing electricity[37]. The electricity into the electrolyser and the hydrogen storage level are both limited by the capacity of electrolyser $C^{el}$, and hydrogen storage $C^S$, respectively:

$$E_{in}^{el}(t) \leq C^{el} \qquad \forall t, \tag{15}$$

$$H^{S_{level}}(t) \leq C^S \qquad \forall t, \tag{16}$$

In addition to the constraints for each component and the flow balances established above, the system can incorporate additional constraints such as temporal correlation and limits on emissions intensity. The temporal correlation constraint requires electricity consumed by the electrolyser to equal dedicated RE generation over certain time intervals, following previous work[14]. In the model, this constraint mandates that the total electricity sold to the grid must be larger than or equal to the total purchased electricity within each time interval:

$$\sum_{t_s}^{t_s+\gamma}(E_{in}^{grid}(t) - E_{out}^{grid}(t)) \geq 0 \qquad \forall t_s \in \theta, \tag{17}$$

where $\gamma$ indicates the length of the time interval required, and $t_s$ marks the beginning of the temporal correlation interval, which falls within the set of time steps $\theta$ tailored to meet the specified time interval. In this work, $\gamma$ ranges from one hour to one year. For example, for a monthly correlation requirement, $\theta$ will include the initial time points of each financial month.

Other policy settings can be incorporated into the model by introducing additional constraints. Additionality is ensured as dedicated RE capacity is installed at the same time as the electrolyser, and geographic correlation is applied by assuming that electrolysers are in the same state as the dedicated RE capacity (i.e. weather data and grid data are state specific). The system can also constrain the emission intensity for hydrogen to below zero to assess the impact on LCOH.

Constraints are also placed on the size of the RE capacity to avoid unphysical results. For example, in states with high electricity prices, the system may grossly oversize the local RE generation capacity to profit from selling electricity to the grid, in essence creating a second business to subsidize hydrogen production. This phenomenon does not align with the practical constraints as the RE generation capacity is usually limited by the available land, project budget and the risk of investment[37]. Thus, we also constrain the $CAPEX$ of the on-grid system such that it does not exceed the $CAPEX$ of the optimized off-grid system without grid-connection.

***Methods for calculating the electricity emissions***

The emissions associated with the electricity used for hydrogen production are calculated using four different methods. Two of these methods are defined by the Australian PGO scheme[20], differentiated by whether they include the buying and selling of RECs, known as the market-based method, or use a location dependent emissions factor defined by the Australian government[24], known as the location-based method. The other two methods take advantage of the ability of our model to assign time and location dependent MEFs or AEFs to the net electricity consumption at each point in time.

For the PGO market-based method, electricity emissions ($e_{GO_M}$) from time $t_1$ to $t_2$ are calculated using the applicable renewable power percentage ($ARPP$), which gives the fraction of the electricity taken from the grid that comes from RE source, and the residual mix factor ($RMF$), which is the estimated average emission factor for the grid that excludes all the electricity generation claimed by RECs[43]. The emissions of electricity used to produce hydrogen is then

$$e_{GO_M} = \left(Q_{grid} \times (1 - ARPP) - Q_{RECs} \times 1000\right) \times RMF. \quad (18a)$$

$$e_{GO_M} = \left(\sum_{t=t_1}^{t_2}(E_{out}^{grid}(t) \times (1 - ARPP) - E_{in}^{grid}(t)\right) \times RMF. \quad (18b)$$

Eq. 18a shows the general case in which the residual electricity, calculated as the total purchased or acquired electricity that is not classified as renewable minus the eligible RECs, is multiplied by the RMF[43], where $Q_{grid}$ is the quantity of grid electricity used and $Q_{RECs}$ is the quantity of RECs surrendered, and Eq. 18b shows how we apply the market-based method in our model from time $t_1$ to $t_2$. We assume that each MWh of RE sold to the grid generates one RECs such that, $Q_{RECs} = E_{in}^{grid}(t)/1000$, which can be surrendered to cover emissions from using grid electricity. If more eligible RECs are surrendered to estimate electricity emissions than the total required to reach zero emissions, $e_{GO_M}$ is equal to zero[24]. Here, we assume that excess RECs will be traded to other facilities, enabling them to claim reductions in their electricity emissions under the market-based method. Consequently, one unit of eligible RECs can represent an emissions reduction on the grid of 0.81 $tCO_2e$ in the market-base method.

To calculate the electricity emissions using the location-based method ($e_{GO_L}$), the location dependent emissions factor for each state ($EF_j$) is defined by the Australian government as part of the National Greenhouse Gas Accounting Factors[24]. The location-based emissions factor is based on the AEFs in each state averaged over the year[43] to calculate the emissions intensity of the electricity used to produce hydrogen:

$$e_{GO_L} = \sum_{t=t_1}^{t_2}(E_{out}^{grid}(t) - E_{in}^{grid}(t)) \times EF_j \quad (19)$$

If the amount of $E_{out}^{grid}(t) < E_{in}^{grid}(t)$, then the system contributes to the reduction of grid GHG emissions.

We define two other emissions calculation methods, using the time and location dependent MEFs or AEFs and the net electricity consumption at each point in time. The MEF-based and AEF-based method keep track of the electricity emissions ($e_{MEF}$ and $e_{AEF}$) from buying and selling electricity to the grid from time $t_1$ to $t_2$ and are defined as

$$e_{MEF} = \sum_{t=t_1}^{t_2}((E_{out}^{grid}(t) - E_{in}^{grid}(t)) \times MEF_j(t)) \quad (20)$$

$$e_{AEF} = \sum_{t=t_1}^{t_2}((E_{out}^{grid}(t) - E_{in}^{grid}(t)) \times AEF_j(t)) \quad (21)$$

Where the MEFs and AEFs of the grid in Australian states $j$ connected to the NEM at each time $t$ are represented as $MEF_j(t)$ and $AEF_j(t)$, respectively. In both cases, the values can be negative, representing the case when more RE is sold to the grid than is bought in each time, and interaction of the system grid with the system is reducing the grid emissions.

## Modelling the Australian National Electricity Market

The NEM is made up of networks in five regions (roughly corresponding to the states) interconnected by a transmission network, including Queensland (QLD), New South Wales (NSW), South Australia (SA), Victoria (VIC) and Tasmania (TAS). Electricity is traded within and between states by the central dispatch engine run by the Australian Energy Market Operator (AEMO)[56]. The spot prices of electricity in each state are determined by the highest bid accepted by the AEMO from a generator to fulfill demand for each 5-minute interval. The historical data on spot prices and emission factors used in our study is from 2023 sourced from AEMO's dataset NEMweb at five-minute intervals[35], and subsequently averaged to one-hour timesteps. Here, a representative historical time series is shown in **Figure 6a**, depicting hourly AEFs and MEFs for Queensland grid profile on 18 Feb 2023. The associated hourly spot prices for the same time and location are shown in **Figure 6b**, along with a bar chart indicating how frequently each type of generator (coal, gas, renewables) becomes the marginal generator over the five-minute intervals in each hourly dataset.

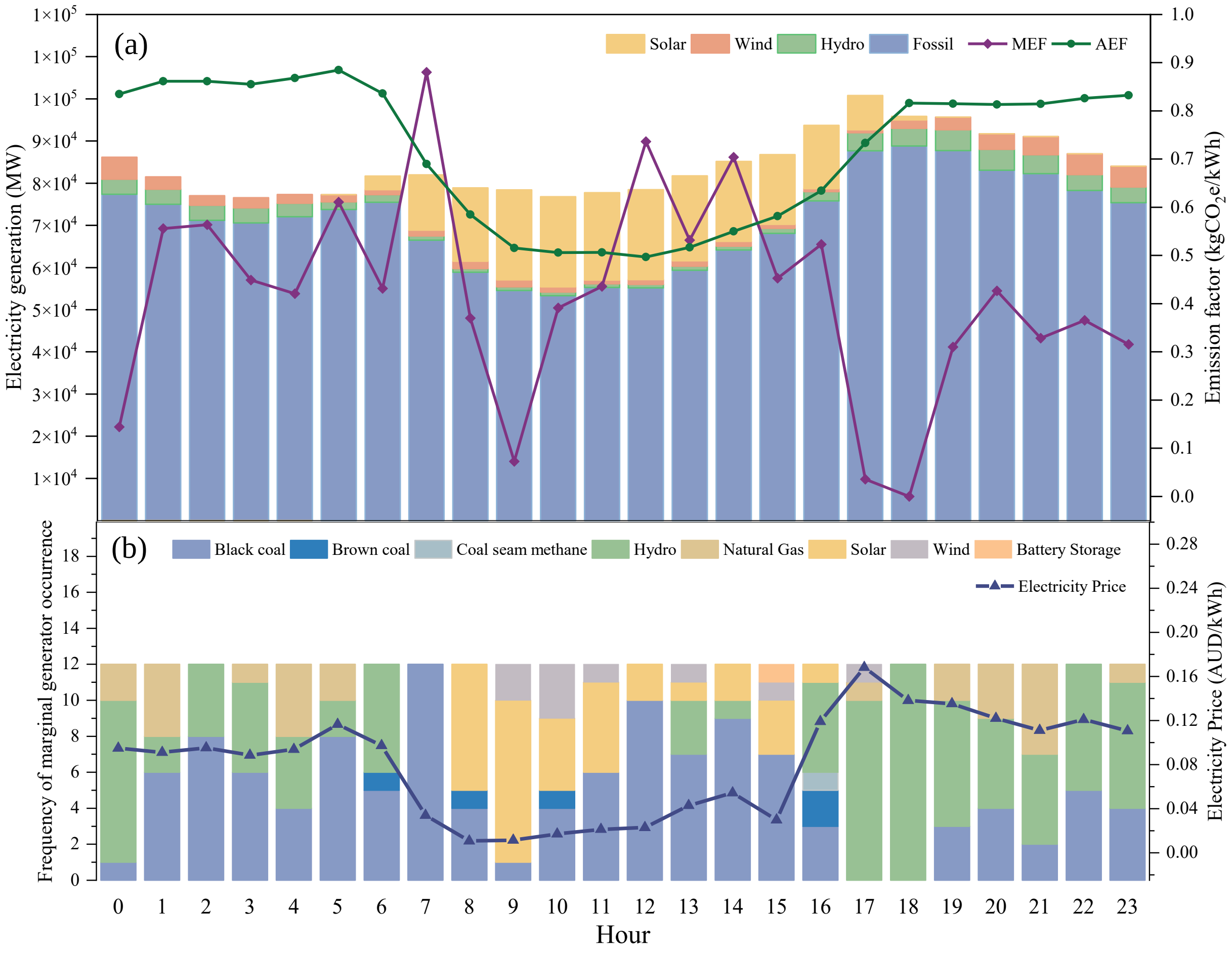

**Figure 6.** The example of the grid profile in QLD on 18 Feb 2023. **(a)** shows marginal emissions factors (MEFs) (purple curve) and average emissions factors (AEFs) (green curve) in relation to the electricity generation of the grid (bars) in each hour. **(b)** shows the frequency of each type of fuel generator becoming marginal (bars), along with electricity prices (blue curve) for each hour.

The graphs show that AEFs reflect the generation mix of the grid, with high AEFs during the early morning and evening peaks due to the reliance on fossil fuel-generated electricity, and lower AEFs during midday when the proportion of RE generation on the grid increases. In contrast, the MEF values are determined by the emissions intensity of the marginal generators and so can change rapidly across the day. The two emission factors can vary widely over some time periods. For example, in hour 18 the AEF is 0.82 kg$CO_2$e/kWh due to the generation mix including RE and fossil sources. At the same time, the MEF plummets to zero, as hydro is the marginal generator and will be used to meet any additional demand. In this work, we consider MEFs to best represent the actual electricity emissions caused by the additional demand on the grid from our system. For example, if our system purchases electricity at hour 18, the generators will be hydro as shown in **Figure 6b**, resulting in no emissions from our purchase. In addition, it tells us which generators would be displaced by RE sold from our system into the grid (and thus quantifies avoided emissions). For instance, if the system sells 1 kWh of renewable electricity at hour 7, the displaced generators will be black coal, reducing grid emissions by 0.88 $kgCO_2e$. Further discussion about the comparison between MEFs and AEFs is given in the **NEM data** section below**.**

It is important to understand how the generator mix and emission profiles of the state grids differ, and how this is reflected in the AEF and MEF values. Yearly average AEFs and MEFs for each state are given in **Table 7,** and the energy generation mix as quantified by the percentage contribution of different generators to annual generation is shown in **Figure 13**. Local grids in QLD, VIC and NSW are dominated by coal power (~60%)[26], with RE meeting a significant portion of the remaining demand, (also given in **Table 7** as 28%, 42%, 31% respectively) leading to average AEFs of 0.69, 0.73, and 0.63 $kgCO_2e/kWh$ in 2023. In contrast, 71% of SA electricity was provided by renewables in 2023, with the rest coming from gas power and imports from other states[26], resulting in a much lower average AEF of 0.23 $kgCO_2e/kWh$. The average AEF in TAS was lower still, at 0.12 $kgCO_2e/kWh$, as over 90% of electricity generation was renewable, mostly from hydro (~70%), followed by wind (~20%).

The average MEF values vary much less between the states, from 0.4-0.52 $kgCO_2e/kWh$, with the exception of TAS, which has a much lower average MEF of 0.19 $kgCO_2e/kWh$. This is because the marginal generators that respond most frequently to increased demand in each state are similar: coal (33-50% of the time), followed by hydro (25-32%), as shown in **Figure 13**. Since TAS has abundant hydro, it is the marginal generator over 75% of the time. The high frequency of hydro as the marginal generator leads to MEFs that are lower than AEF in heavily coal-reliant states (QLD, NSW, VIC), while the fossil fuel technologies are marginal generators 45% of the time in SA, resulting in larger MEF than AEF.

## Modelling geographic correlation

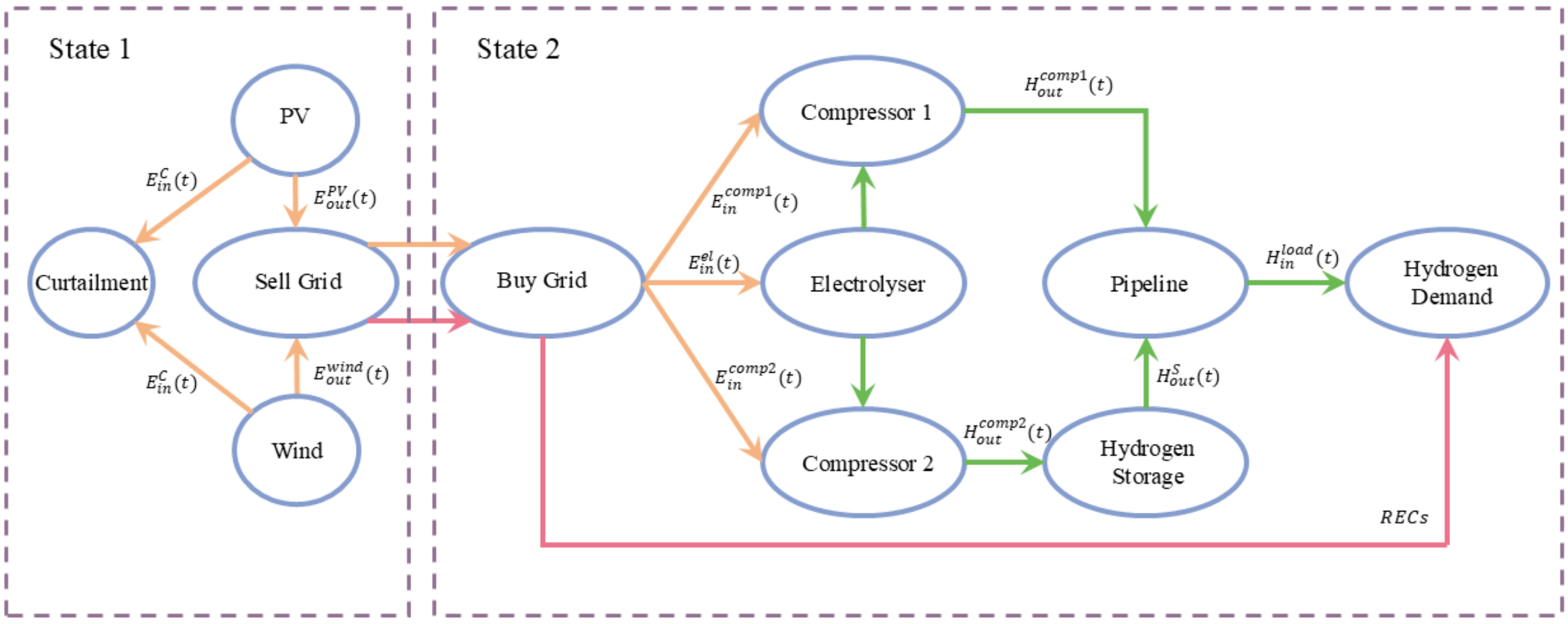

**Figure 7** Overview of hybrid grid-connected hydrogen production optimization model without the geographic correlation constraints. Ovals represent nodes, which are system components responsible for transformation of energy commodities. Lines represent edges which indicate the flow of an energy commodity and arrows represent the flow directions. Orange lines indicate electricity flow (kW), while green lines represent hydrogen flow (kg/hour), and pink lines show renewable energy certificates (RECs) associated with renewable electricity integrated into the grid. The dashed lines indicate the boundary between two states. The Sell Grid node locates in State 1, where the renewable generators (wind farm and PV field) integrate electricity to the grid, while the hydrogen production system and the Buy Grid node are situated in State 2. State 1 and State 2 are interconnected two bidding zones via a transmission line.

Apart from the scenario in which renewable generators (the wind farm and solar field) are co-located with the electrolyser, meeting the geographic correlation requirement as shown in **Figure 5**, our study uses the same model but with modified constraints to study a scenario where the electrolyser is not co-located with renewable generators (without the direct line connection) to investigate the policy impact of geographic correlation on the emissions and cost of hydrogen production as shown in **Figure 7**. Under this scenario, we assume that the hydrogen producer invests in self-owned renewable generators located in states different from where the hydrogen production system is situated[57]. All renewable electricity generated will be sold to the regional grid (the Sell Grid node), generating profit and producing RECs in the process. The hydrogen production system is grid-connected with the regional grid (the Buy Grid node), leveraging the NEM's transmission network to match its electricity consumption with the self-invested renewable generators using RECs information.

Distinct from Eq. (14) above, the electricity cost over one year includes costs associated with buying electricity from the Buy Grid and negative "costs" from selling electricity to the Sell Grid in each time $t$ and is given by,

$$C^e = \sum_{t=0}^{8759}(E_{out}^{grid}(t) \times (P_{buy}(t) + TS) - E_{in}^{grid}(t) \times P_{sell}(t)), \tag{22}$$

where $P_{buy}(t)$ and $P_{sell}(t)$ are the electricity spot price in the Buy Grid and Sell Grid at each time $t$, respectively. For example, the hydrogen producer establishes a hydrogen production system in Newcastle (NSW), including a grid-connected electrolyser, compressors, hydrogen storage, and transportation pipeline facilities, while developing a wind farm in Burnie (TAS). The electricity purchase cost will be determined by the Buy Grid Electricity price in NSW ($P_{buy}(t)$), while the revenue from electricity sales will be based on the Sold Grid Electricity price in TAS ($P_{sell}(t)$). The RECs generated by the Burnie wind farm will be used to certify the electricity consumed by the hydrogen production system in Newcastle.

## Data input

## Selected locations

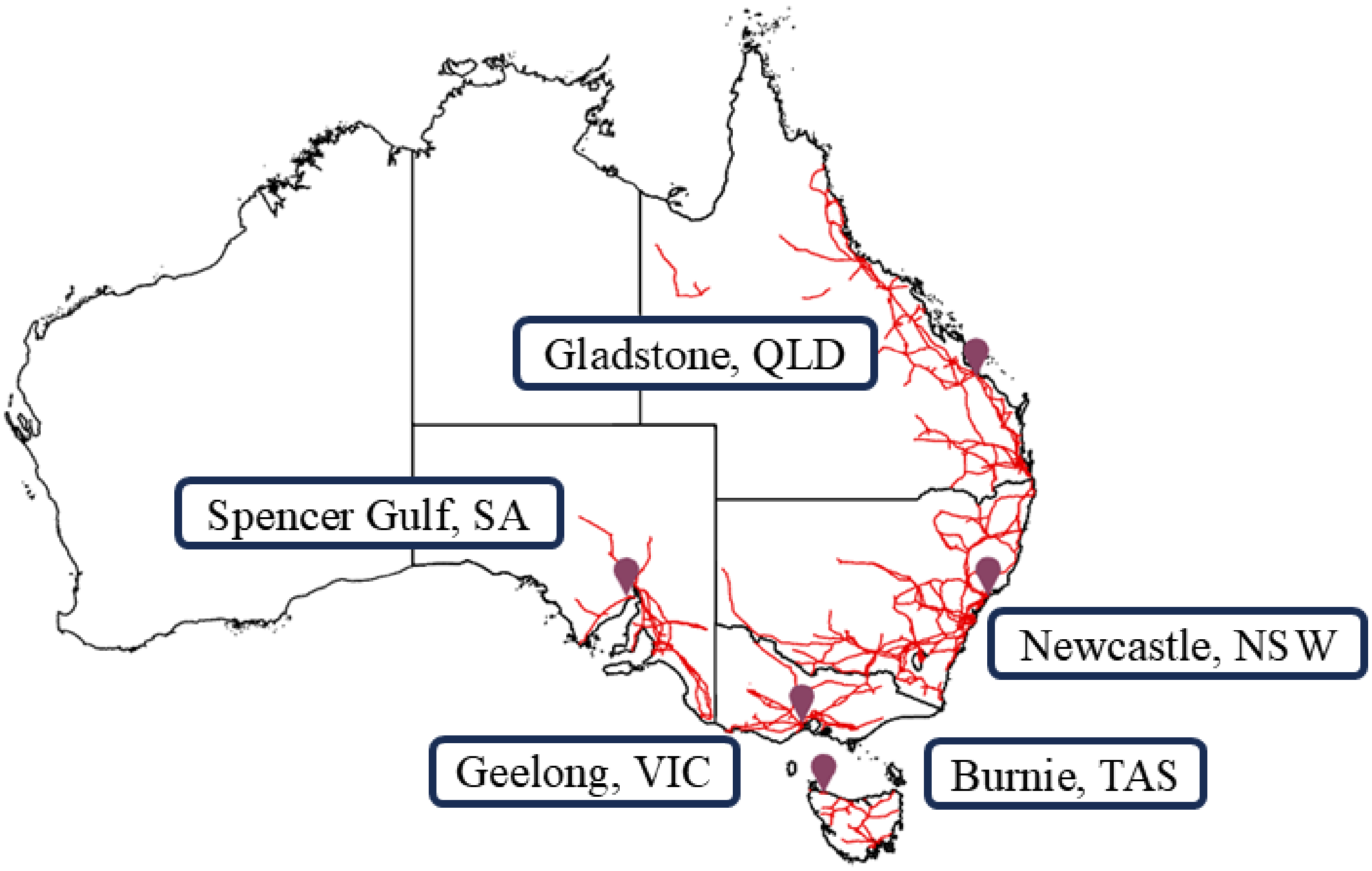

**Figure 8** The five locations selected as the case study.

**Table 2** The latitude and longitude of the five locations

| State | Location | Latitude | Longitude |
|---|---|---|---|
| Queensland | Gladstone | -24.08 | 151.28 |
| South Australia | Spencer Gulf | -32.94 | 137.45 |
| Tasmania | Burnie | -40.89 | 145.25 |
| Victoria | Geelong | -38.07 | 144.39 |
| New South Wales | Newcastle | -32.92 | 151.75 |

In this study, we select one location for a hydrogen production plant in each of the five states participating in the Australian National Electricity Market (NEM) as the example projects using Australian hydrogen projects dataset as the guidance[49]. The database provides a comprehensive overview of all hydrogen projects across Australia, including detailed information such as the energy sources used for hydrogen production, the production methods employed, and the annual hydrogen production amounts for each project. We focus on projects that source their energy from the grid and select those that are also participating in the current trial projects under the PGO scheme[58].

The geographical locations of the five cases we study, all of which have available grid connection conditions, are shown in **Figure 8**, where the purple dots represent the locations of the hydrogen production systems, and the red line indicates the grid electricity transmission line sourced from Geoscience Australia dataset[59]. The specific latitude and longitude are listed in **Table 2**.

## Model input technical and economic parameters

The input technical and economic parameters for the optimization model are chosen from best available data sources in the literature and listed in **Table 3**. In particular, we use data from the latest Gencost report as extensively as possible. This economic report estimates the costs of building new electricity generation and hydrogen production facilities specific to Australia[60]. Costs are given for 2023 and converted to US dollars with an exchange rate of 0.7 US dollars per Australian dollar which is in line with the exchange rate adopted in latest Gencost report[60].

**Table 3** Model input parameters

| Parameters | Expression | Units | Current Value |
|---|---|---|---|
| Capex of electrolyser | $I^{el}$ | USD/kW | 1343.3[60] |
| Capax of wind | $I^{wind}$ | USD/kW | 2126.6[60] |
| Capex of PV | $I^{PV}$ | USD/kW | 1068.2[60] |
| FOM of electrolyser | $FOM^{el}$ | USD/kW | 37.4[61] |
| FOM of wind | $FOM^{wind}$ | USD/kW | 17.5[60] |
| FOM of PV | $FOM^{PV}$ | USD/kW | 11.9[60] |
| VOM of electrolyser | $VOM^{el}$ | USD/kgH$_2$ | 0.02[62] |
| Electricity consumption of pumping unit hydrogen into pipeline | $\mu_{out}^{comp1}$ | kWh/kgH$_2$ | 0.83 |
| Efficiency of electrolyser | $\eta$ | | 70%[53] |
| Higher heating value of hydrogen | $HHV$ | kWh/kg H$_2$ | 39.4[41] |
| Discount rate | $i$ | | 6% |
| Project lifetime | $n$ | | 25 |
| Transmission use of system fee | $TS$ | USD/kWh | 0.007[63] |

In this study, we explore two options of hydrogen storage equipment: pipeline storage and the Lined Rock Cavern (LRC) underground storage. The relationship between the total capex for both the pipeline and LRC underground storage systems and their storage capacity are calculated using the following Eqs. (23)-(24)[64]:

$$log_{10}(I_{pipe}^{S}) = -0.0285 log_{10}\left(\frac{C_{pipe}^{S}}{1000}\right) + 2.7853 \quad (23)$$

$$log_{10}(I_{LRC}^{S}) = 0.217956\left(log_{10}\left(\frac{C_{LRC}^{S}}{1000}\right)\right)^{2} - 1.575209 log_{10}\left(\frac{C_{LRC}^{S}}{1000}\right) + 4.463930 \quad (24)$$

Where unit cost of pipeline storage $I_{pipe}^{S}$ and LRC underground storage $I_{LRC}^{S}$ are in USD per kg of $H_2$, and their capacity is in kg of $H_2$.

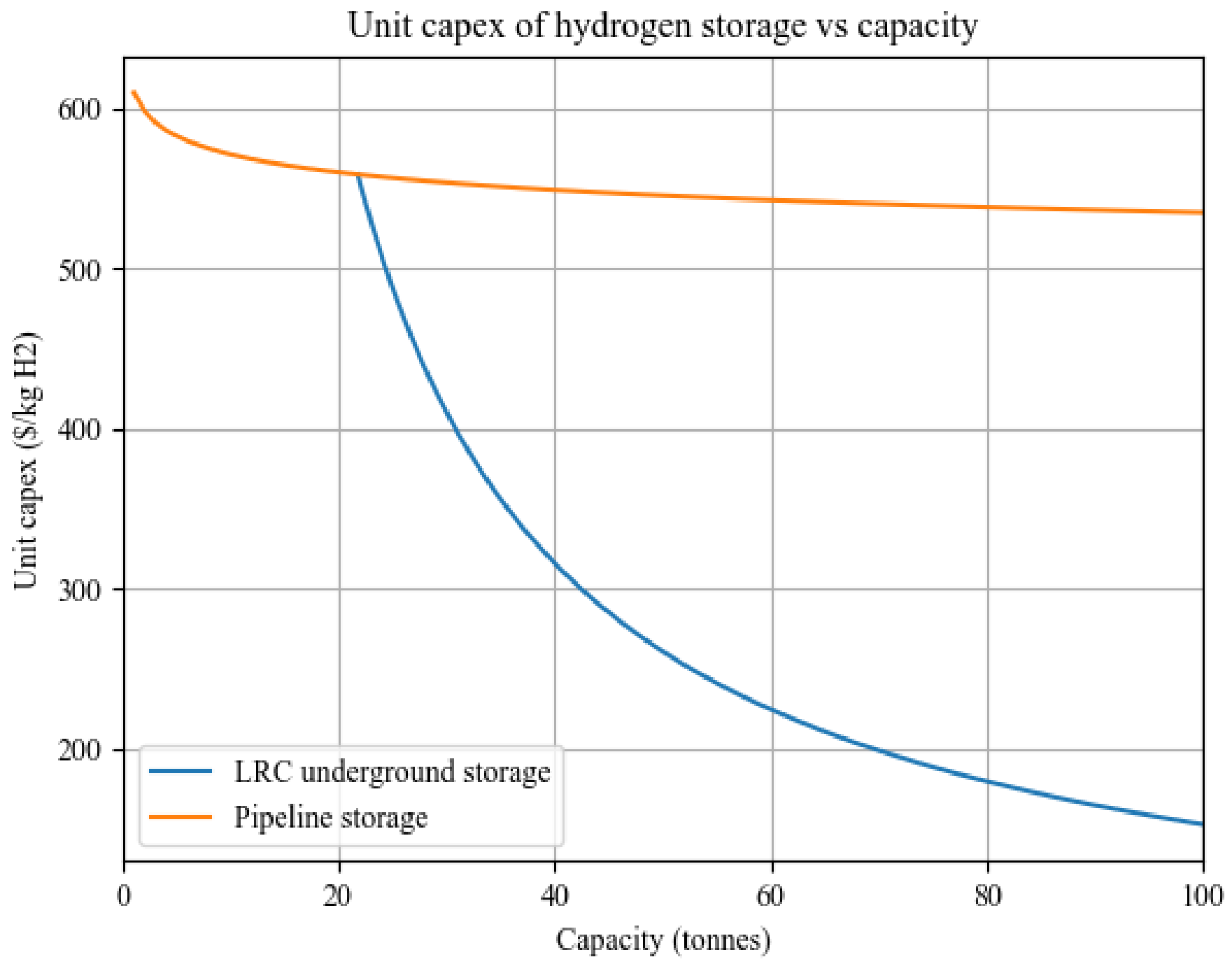

**Figure 9** The unit capex of hydrogen storage as a function of the capacity.

Here, we utilize the piecewise function to linearize these two cost functions to reduce the solving complexity, and the curve of unit capex of those two types of hydrogen storage based on their capacity is shown in **Figure 9.** Noticeably, LRC underground storage capex decreases significantly and become more economical than pipeline storage for larger storage sizes. Hence, we assume the system will opt for pipeline storage if the hydrogen storage size is under 21742kg and the electricity consumption of pumping unit hydrogen into the storage $\mu_{out}^{comp2}$ is 0.83 kWh/kg$H_2$. Otherwise, the system will choose LRC underground storage and $\mu_{out}^{comp2}$ is 1.24 kWh/kg$H_2$, due to the higher-pressure requirement of underground storage (150 bar) compared to the pipeline (100 bar).

## Renewable generation

**Table 4** PV field model parameters.

| Parameter | Value |
|---|---|
| Rated Capacity | 1MW |
| Array type | 2-axis tracking |
| Tilt | 0° |
| Azimuth | 180° |
| System losses | 15% |

**Table 5** Wind farm model parameters.

| Category | Parameter | Value |
|---|---|---|
| Turbine | Model | 2020 ATB NREL Ref. 4 MW |
| | Hub height | 150m |
| | Number of turbines | 80 |
| | Rated Capacity | 320MW |
| Availability losses | Turbine | 3.58% |
| | Balance of plant | 0.5% |
| | Grid | 0.5% |
| Electric losses | Efficiency | 3% |
| Turbine performance losses | Sub-optimal performance | 1.1% |
| | Generic power curve adj. | 1.7% |
| | Site-specific power curve adj. | 0.81% |
| | High wind hysteresis | 0.4% |
| Environmental losses | Degradation | 0.1% |
| | Environmental | 0.4% |

For PV generation, we import the PVWatt module from the NREL's PySAM package[51] to calculate the reference electricity generation of a 1MW PV field for the selected locations. The system's parameters are provided in **Table 4**. For wind power generation, we select the NREL's 2020 ATB Reference 4MW wind turbine with the system parameters provided in **Table 5**. The reference electricity generation of a 320MW wind farm is calculated using NREL PySAM's windpower module[51].

**Table 6** RE capacity factor of five locations in 2021 and 2023.

| Location | Year | RE capacity factor |
|---|---|---|
| Gladstone (QLD) | 2021 | 0.28 |
| | 2023 | 0.29 |
| Spencer Gulf (SA) | 2021 | 0.35 |
| | 2023 | 0.34 |
| Burnie (TAS) | 2021 | 0.55 |
| | 2023 | 0.54 |
| Geelong (VIC) | 2021 | 0.42 |
| | 2023 | 0.39 |
| Newcastle (NSW) | 2021 | 0.24 |
| | 2023 | 0.23 |

**Table 6** provides RE capacity factor of five locations for the years 2021 and 2023. The RE capacity factor ($CF_{RE}$) of each location is calculated using the optimized off-grid system's configuration to evaluate the local RE resources as shown in the following equation:

$$CF_{RE} = \frac{\sum_{t=0}^{8759}\left(E_{out}^{wind}(t)+E_{out}^{PV}(t)\right)}{\left(C^{wind}+C^{PV}\right)\times 8760} \quad (25)$$

Where $E_{out}^{wind}(t)$ and $E_{out}^{PV}(t)$ are the wind and PV generation at each time $t$, and $C^{wind}$ and $C^{PV}$ represent the capacity of wind and PV, respectively.

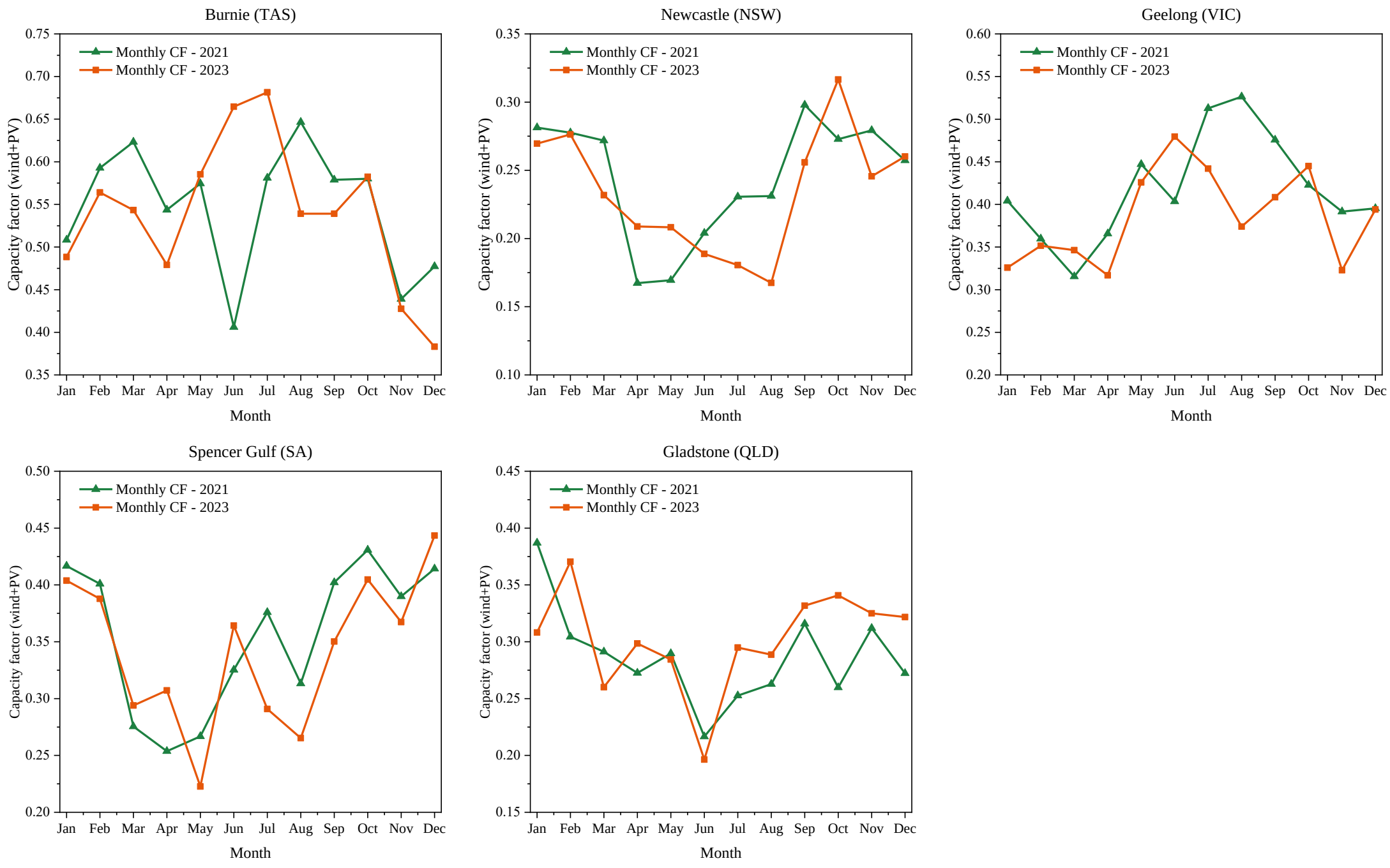

**Figure 10** Monthly average capacity factor in five locations in 2021 and 2023.

**Figure 10** visualizes the monthly average RE capacity factors for the five cases studied across five states in the NEM, based on the optimized off-grid system configuration. The different coloured curves represent the monthly average RE capacity factors for the year 2021 and 2023.

## NEM data

**Table 7** Regional grid information of five locations in 2021 and 2023.

| Location | Year | Regional yearly average electricity prices (AUD/MWh) | Regional grid RE penetration | Yearly average MEF of regional grid ($kgCO_2e/kWh$) | Yearly average AEF of regional grid ($kgCO_2e/kWh$) |
|---|---|---|---|---|---|
| Gladstone (QLD) | 2021 | 89 | 19% | 0.63 | 0.76 |
| | 2023 | 95 | 28% | 0.47 | 0.69 |
| Spencer Gulf (SA) | 2021 | 59 | 63% | 0.51 | 0.3 |
| | 2023 | 92 | 71% | 0.4 | 0.23 |
| Burnie (TAS) | 2021 | 35 | 102% | 0.26 | 0.09 |
| | 2023 | 57 | 93% | 0.19 | 0.12 |
| Geelong (VIC) | 2021 | 48 | 34% | 0.54 | 0.8 |
| | 2023 | 64 | 42% | 0.43 | 0.73 |
| Newcastle (NSW) | 2021 | 73 | 23% | 0.61 | 0.74 |
| | 2023 | 98 | 31% | 0.52 | 0.63 |

**Table 7** provides regional grid profiles of five locations for the years 2021 and 2023. The data of RE penetration in the regional grid are sourced from Open Electricity's Tracker[26] and yearly electricity prices, MEFs, and AEFs for the regional grid are sourced from the NEMweb dataset, which provides data at 5-minute intervals, and are averaged across the year[35].

It is worth noting that the grid profile during 2021 and 2023 for each of the five states reveals significant electricity prices variations across different years. For instance, in 2023, SA yearly average electricity price surged to 92 AUD/MWh, much higher compared to 2021 (59 AUD/MWh). The yearly average MEFs in five states are all decreasing, and similar decreasing trend can be seen in the yearly average AEFs due to with higher grid RE penetration in NSW, VIC, SA and QLD except for TAS, where RE penetration is reducing from 102% in 2021 to 93% in 2023 and the yearly average AEFs increase from 0.09 $kgCO_2e/kWh$ in 2021 to 0.12 $kgCO_2e/kWh$ in 2023.

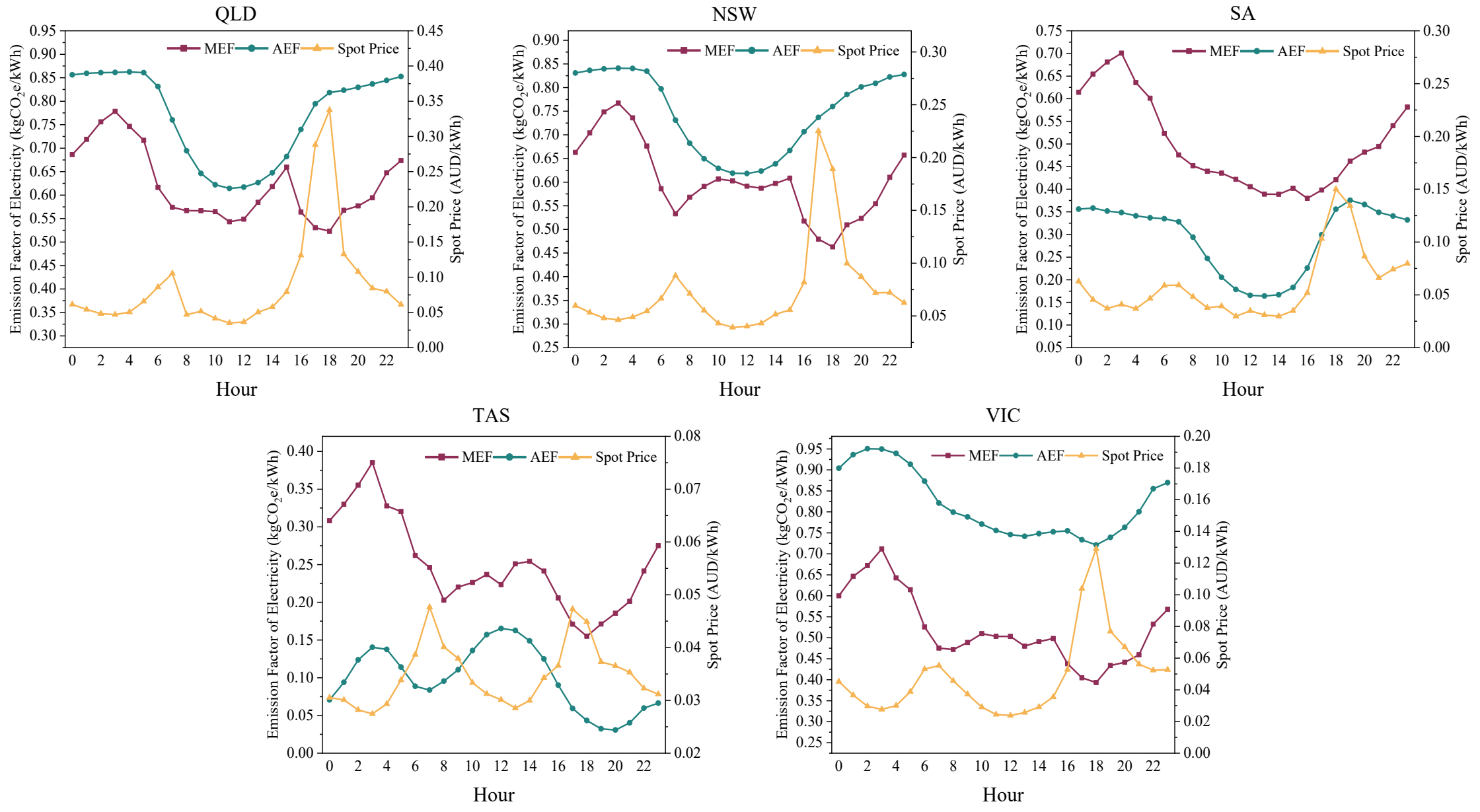

**Figure 11** Hourly average of MEFs vs AEFs and spot prices of the regional grid in five states in 2021.

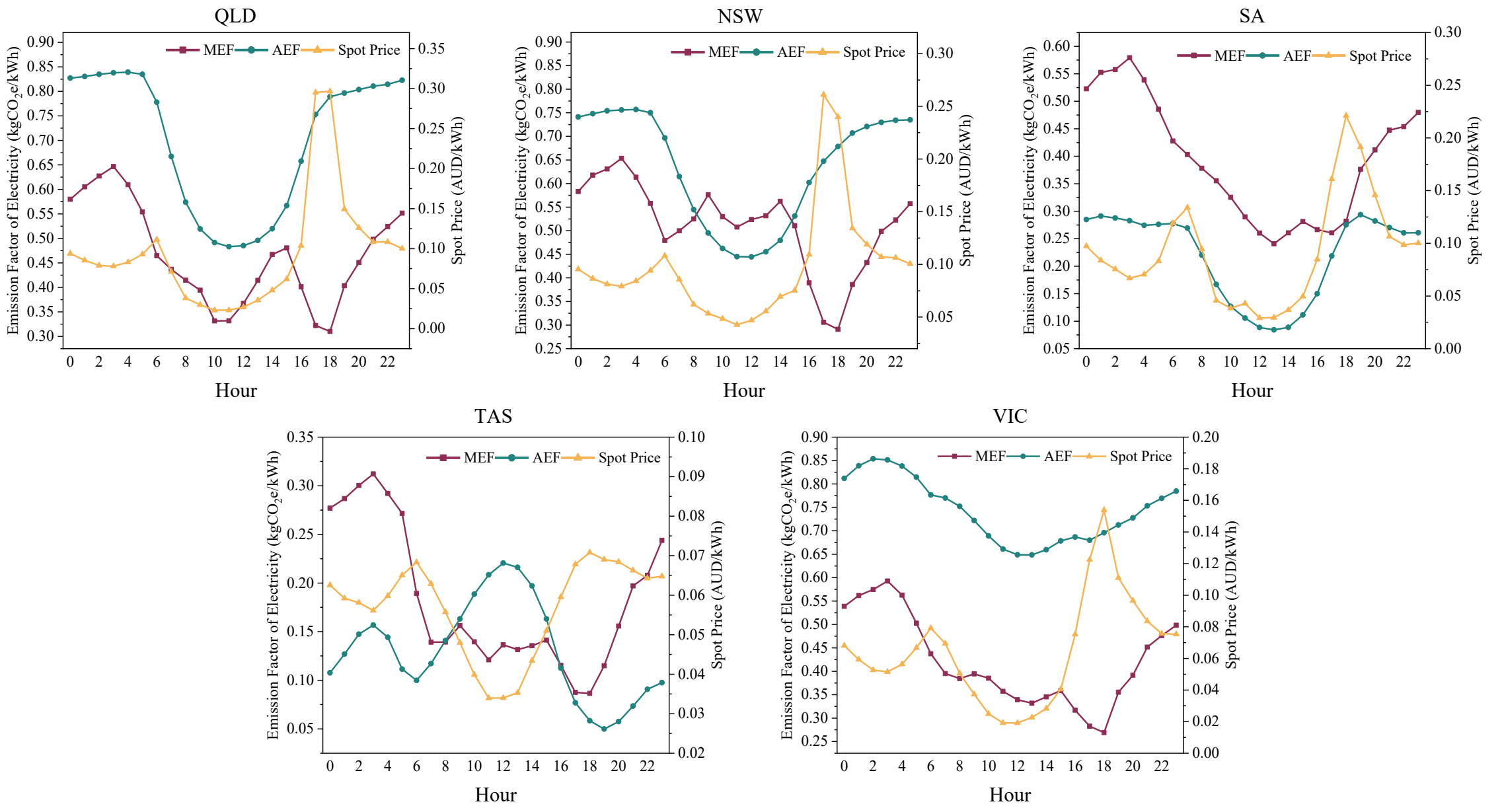

**Figure 12** Hourly average of MEFs vs AEFs and spot prices of the regional grid in five states in 2023.

We simulate the interaction between the hydrogen production system and the grid via modelling the grid using the historical electricity price, MEF and AEF data published by Australian Energy Market Operator (AEMO). With the help of the python package NEMOSIS[65], we get access to hourly electricity prices of five regions in the NEM. We obtain MEF utilising python package NEMED[66] by identifying the marginally dispatched generators from AEMO's price setter files and associating each marginal generator with its respective carbon dioxide equivalent intensity index published in AEMO's NEMweb dataset for each time interval[35], AEF is computed by the energy sent out consumed from the grid (including the inter-region electricity trade) and its associated emissions. Emissions data is obtained using NEMED, and data of energy sent out and the inter-region electricity trade is acquired via NEMOSIS.

Through comparing the hourly average of MEFs vs AEFs in 2021 and 2023 as shown in **Figure 11** and **Figure 12**, it is noted that the hourly average MEFs are quite different with AEFs. In this work, we adopt the two methods: AEF-based method and MEF-based method, to compare the electricity emissions calculated by two PGO methodologies. The two emissions factors, AEFs and MEFs, answer distinct questions. If we assume that the electricity demand in question represents a portion of the existing demand, then AEFs are more appropriate for answering the question of what share of grid emissions this electricity demand is responsible for. Conversely, if the electricity demand under analysis represents an incremental increase or decrease, MEFs are preferred to determine the impact of the newly changed demand on grid emissions[38]. This is because changes in demand are not distributed proportionally across all generators within the electricity system. In real operation, specific generators respond to these electricity demand change, and it is the emission intensity of these generators (MEFs) that drives the resulting emission variations[39]. Thus, we consider the MEF-based method to best represent actual emissions in this work because we use historical grid data to assess emission variations arising from the additional electricity demand of the hybrid grid-connected hydrogen production systems and study how policy interventions like temporal and geographic correlations requirements may impact emissions.

Previous research has emphasised that MEFs are more accurate than AEFs for calculating emissions resulting from shifts in electricity demand and for studying the effects of policy interventions. For example, Sengupta et al. (2022) pointed out MEFs provide accurate assessments of emissions, particularly when analysts or policymakers evaluate interventions that result in changes to electricity demand[67]. Holland et al. (2022) discovered that using average rather than marginal emissions to predict the impacts significantly overestimates the emission benefits of the policy intervention of Biden administration's 2030 EV target. They also highlighted it is estimates of marginal emissions that are needed to accurately evaluate the impacts of policies or behaviours that cause changes in the demand or supply of electricity[68]. Elenes et al. (2022) found that average emission factors have lower accuracy when estimating emissions from demand shifts. In contrast, incremental marginal emission factors (emissions in response to a small shift in electricity demand) can reproduce the emission changes of a power grid model under many testing conditions and scenarios[69]. Li et al. (2017) indicated that for evaluating the emission implications of policy and technology interventions, such as electric vehicle tax credits and energy storage solutions, MEFs are a more suitable metric than AEFs[70]. Siler-Evans et al. (2012) revealed that AEFs can substantially misestimate the avoided emissions resulting from an intervention, in contrast to MEFs[71].

**Table 8** Scope 2 emission factors and residual mix factors for consumption of electricity and applicable renewable power percentage.

| State | Scope 2 emission factor ($kgCO_2e/kWh$) | Residual mix factor ($kgCO_2e/kWh$) | ARPP |
|---|---|---|---|
| Queensland | 0.71 | 0.81 | 18.72% |
| South Australia | 0.23 | 0.81 | 18.72% |
| Tasmania | 0.15 | 0.81 | 18.72% |
| Victoria | 0.77 | 0.81 | 18.72% |
| New South Wales | 0.66 | 0.81 | 18.72% |

**Table 8** provides the scope 2 emission factors of five states used in the PGO location-based method, and the residual mix factor and applicable renewable power percentage applied in the PGO market-based method, sourced from the latest National Greenhouse and Energy Reporting (Measurement)

Determination 2008[24]. The scope 2 emission factors are state-based, calculated as the average emission factor for all electricity consumed from the grid within a specific state, territory, or electricity grid. The applicable renewable power percentage ($ARPP$) is calculated by[72],

$$ARPP = (E^{RE} + E^{ad})/(E^{AC} - E^{ex}) \quad (26)$$

Where $E^{RE}$ is renewable electricity required for the year, $E^{ad}$ is the cumulative adjustment, $E^{AC}$ is the estimated relevant acquisitions of electricity for the year and $E^{ex}$ represents estimated exemptions. Following the example provided in Australian National Greenhouse Accounts Factors 2024[43], we adopt the average $ARPP$ value for the calendar years 2023 and 2024 in this work[72], which is 18.72%. Another key parameter used in the PGO market-based method is the residual mix factor ($RMF$), which is calculated nationally by adjusting the national location-based scope 2 emission factor to exclude the emissions benefits of all claimable renewable generation through LGCs[43], calculated by,

$$RMF = E / (Q - (RECs \times 1000)) \quad (27)$$

Where $RMF$ is the residual mix factor for the financial year. $E$ represents the emissions from on-grid generation based on NGER reports for the preceding reporting year. $Q$ is the electricity generation sent out to the grid. $RECs$ represent renewable energy certificates created for renewable energy generation which occurred in the same financial year. The spot price of RECs is estimated to range from 20 AUD to 60 AUD per unit, based on spot price variations observed during 2023-2024[33].

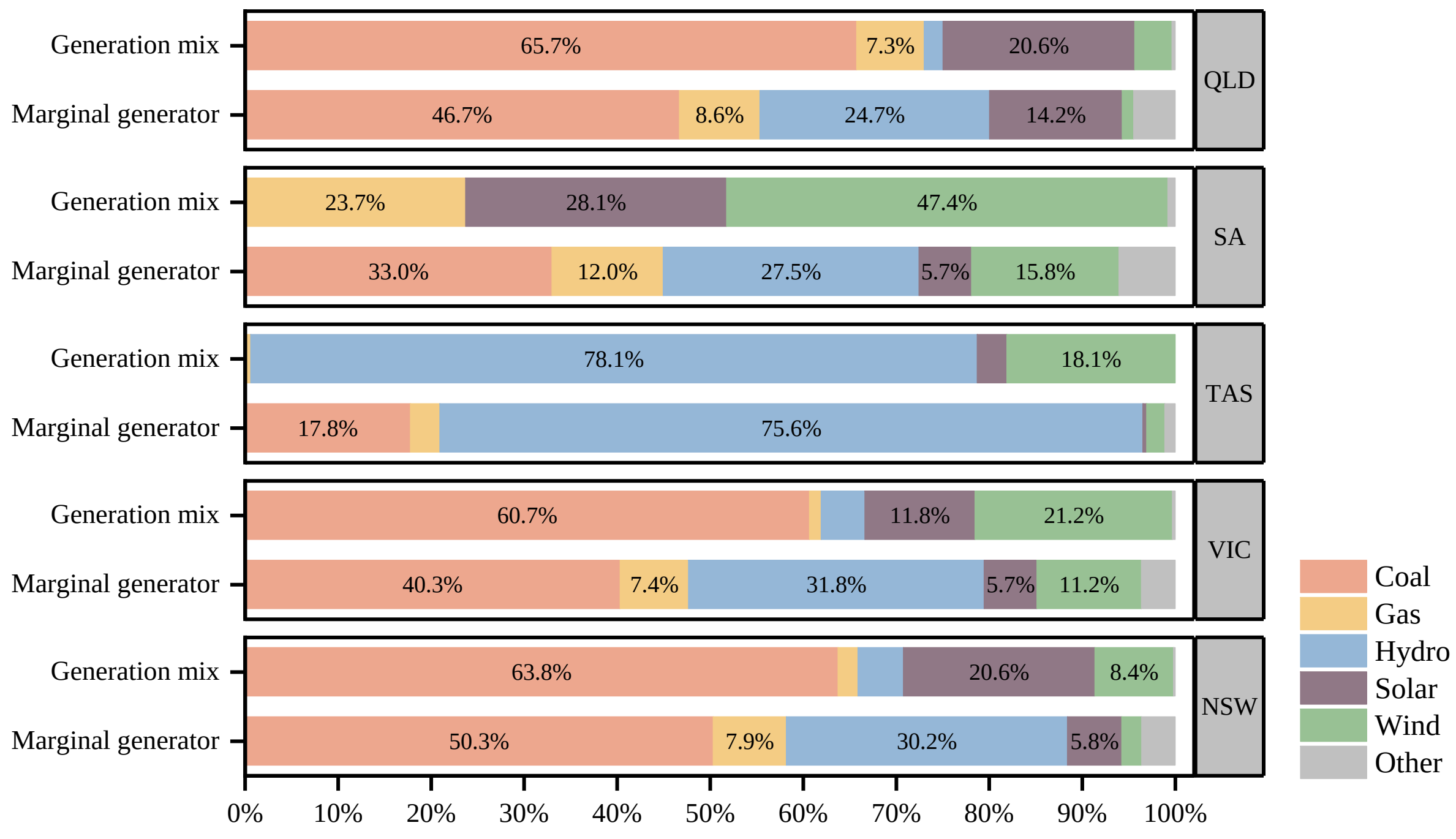

**Figure 13** Electricity generation mix and frequency of marginal generator fuel types across 5-min intervals of five states in 2023

## Time series illustrating the modelled system operation

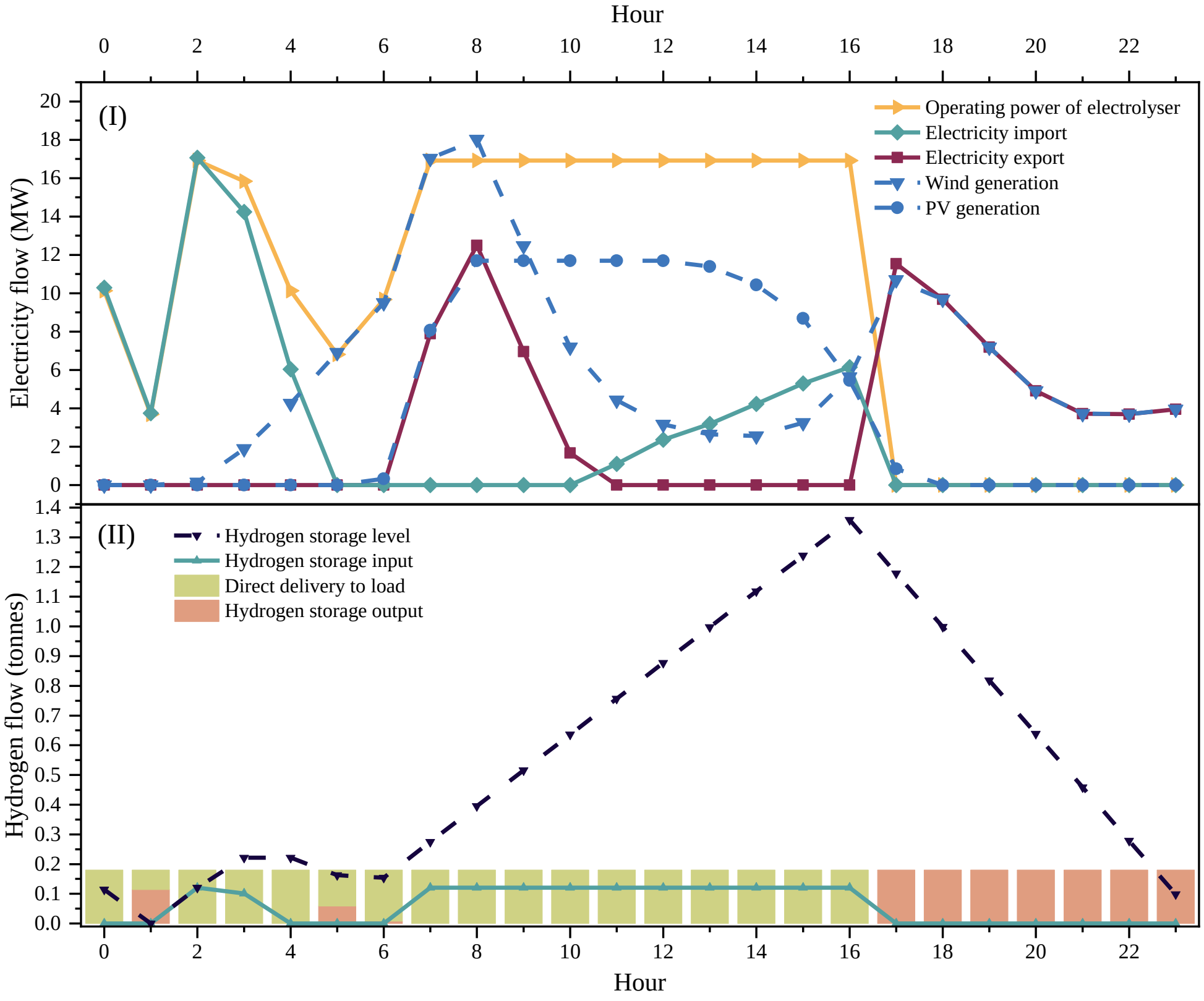

**Figure 14** Example of time resolved output of the model over 24 hours in the case of Gladstone (QLD) on 9 Sep 2023. **(I)** Electricity flow [MW] from different components and **(II)** hydrogen flow [in tonnes].

**Figure 14** gives a representative time series illustrating the modelled system operation on 9 Sep 2023 in Gladstone (QLD). By tracing the electrical flow from various components in **Figure 14I** and the hydrogen flow in of **Figure 14II**, we can monitor the operation of the electrolyser, whether power is provided by local RE or the grid, and the hydrogen storage level.

At time 0 (12 am), local RE generation is minimal and the hydrogen storage does not have sufficient inventory to meet the demand. Thus, the system imports electricity from the grid to power the electrolyser, producing hydrogen to meet demand. Once the local RE generation is increasing from 3 am to 5 am, the system gradually reduces its import of grid electricity, instead using local RE to operate the electrolyser. From 7 am to 10 am, as wind and PV generation rises, the system not only utilises local RE to produce hydrogen for meeting demand and replenishing hydrogen storage but also exports excess electricity to the grid for profit. The system resumes importing the electricity at 11 am to maintain the electrolyser at nominal power as wind generation decreases and continues importing until 16 pm due to low electricity prices during this period. Then, during the high electricity prices of the day in the evening from 6 pm to 11 pm, the system chooses to cease hydrogen production from the electrolyser and export all locally generated electricity to the grid. At this period, the hydrogen stored in hydrogen storage is transported to meet demand.

### Limitations

Firstly, our model optimizes the system configuration over a single year with perfect foresight, and we do not investigate how natural variation in weather patterns and the dynamic changes in grid spot price and emissions over the years would affect optimal sizing or operation. Secondly, we use fixed estimates for component costs and do not include location dependent cost factors. Thirdly, we do not consider additional revenue stream for the producer from selling the RECs, due to the high uncertainty of the future fluctuations of the RECs prices[73]. Instead, we focus on the comparative emissions and costs of systems in different states, under different constraints to assess the accuracy of PGO certification methodologies and the efficacy of different policy settings in minimizing emissions while maintaining low costs.

# Additional Results

### AEF-based method vs Location-based method

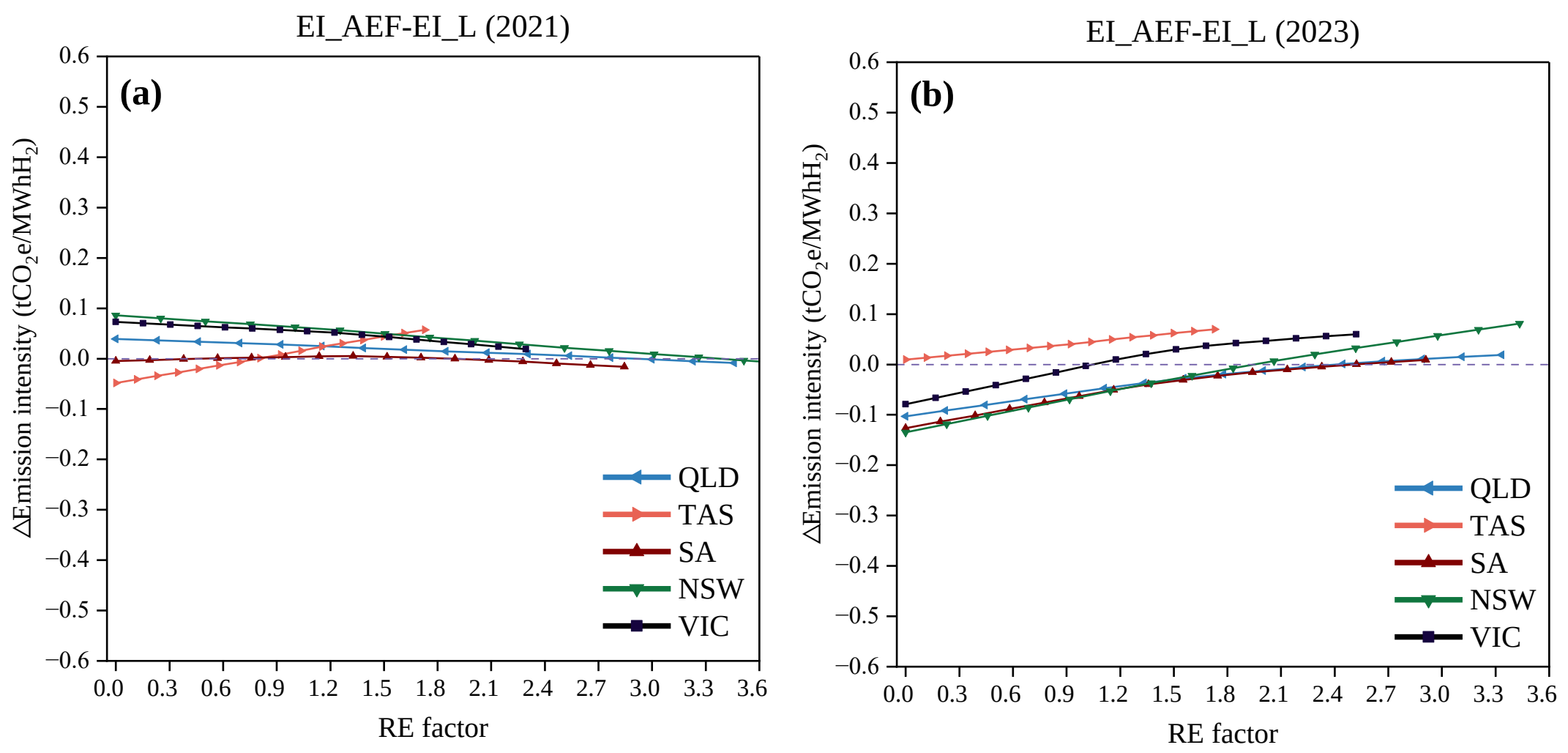

**Figure 15 (a)-(b)** The difference in emission intensities calculated by AEF-based method and PGO location-based method for different RE factors across the five states in 2021-2023.

**Figure 15 (a)-(b)** shows the differences in emission intensities during 2021 and 2023, calculated by subtracting the PGO location-based method from the AEF-based method. The 2021 results are included here to validate the findings. A good alignment between the PGO location-based method and the AEF-based method can be observed across all systems investigated during 2021 and 2023, with no distinct trends. The differences in emission intensities (EI_AEF - EI_L) remain within a narrow range of -0.15 to 0.1 $tCO_2e/MWh_{H2}$. This demonstrates that PGO location-based method accurately reflects the average GHG emissions associated with electricity consumption from a specific regional grid, and the temporal resolution of the average emission factors does not significantly affect the emissions estimates calculated by this approach.

**Table 9** Selected numerical results of **Figure 2**

| Location | RE factor | EI_MEF ($kgCO_2e/kgH_2$) | EI_Market ($kgCO_2e/kgH_2$) | EI_RECs ($kgCO_2e/kgH_2$) | EI_L ($kgCO_2e/kgH_2$) | $D_L$ | $D_{Market}$ |
|---|---|---|---|---|---|---|---|
| Gladstone | 0 | 27.61 | 37.6 | 0 | 40.55 | -47% | -36% |
| (QLD) | 3.33 | -13.62 | 0 | -26.02 | -21.41 | -57% | 91% |
| Spencer | 0 | 22.55 | 37.6 | 0 | 13.14 | 42% | -67% |
| Gulf (SA) | 2.91 | -9.15 | 0 | -20.28 | -5.36 | 41% | 122% |
| Burnie | 0 | 10.06 | 37.6 | 0 | 8.57 | 15% | -274% |
| (TAS) | 1.73 | -6.54 | 0 | -22.57 | -3.82 | 42% | 245% |
| Geelong | 0 | 26.9 | 37.6 | 0 | 43.98 | -63% | -40% |
| (VIC) | 2.52 | -9.65 | 0 | -20 | -17.54 | -82% | 107% |
| Newcastle | 0 | 32.42 | 37.6 | 0 | 37.7 | -16% | -16% |
| (NSW) | 3.43 | -14.07 | 0 | -27.03 | -20.59 | -46% | 92% |

The conversion between ($kgCO_2e/kgH_2$ and $tCO2e/MWh_{H2}$ is based on the higher heating value of hydrogen (39.4 $kWh/kgH_2$[41]). The difference percentage between EI_MEF and EI_Market ($D_{Market}$), between EI_MEF and EI_L ($D_L$) is calculated by:

$$D_{Market} = \frac{(EI_MEF - EI_Market)}{|EI_MEF|} \tag{28a}$$

$$D_{Market} = \frac{(EI_MEF - EI_RECs)}{|EI_MEF|} \tag{28b}$$

$$D_L = \frac{(EI_MEF - EI_L)}{|EI_MEF|} \tag{29}$$

When the system generates more RECs that result in the emission intensity calculated by the market-based method falling below zero, EI_Market is assumed to be zero[24]. Then we use the emission offset by RECs (EI_RECs) as the emission intensity of hydrogen calculated by the market-based method to compare the normalised difference scale with the EI_MEF ($D_{Market}$).

## MEF-based method vs PGO methods

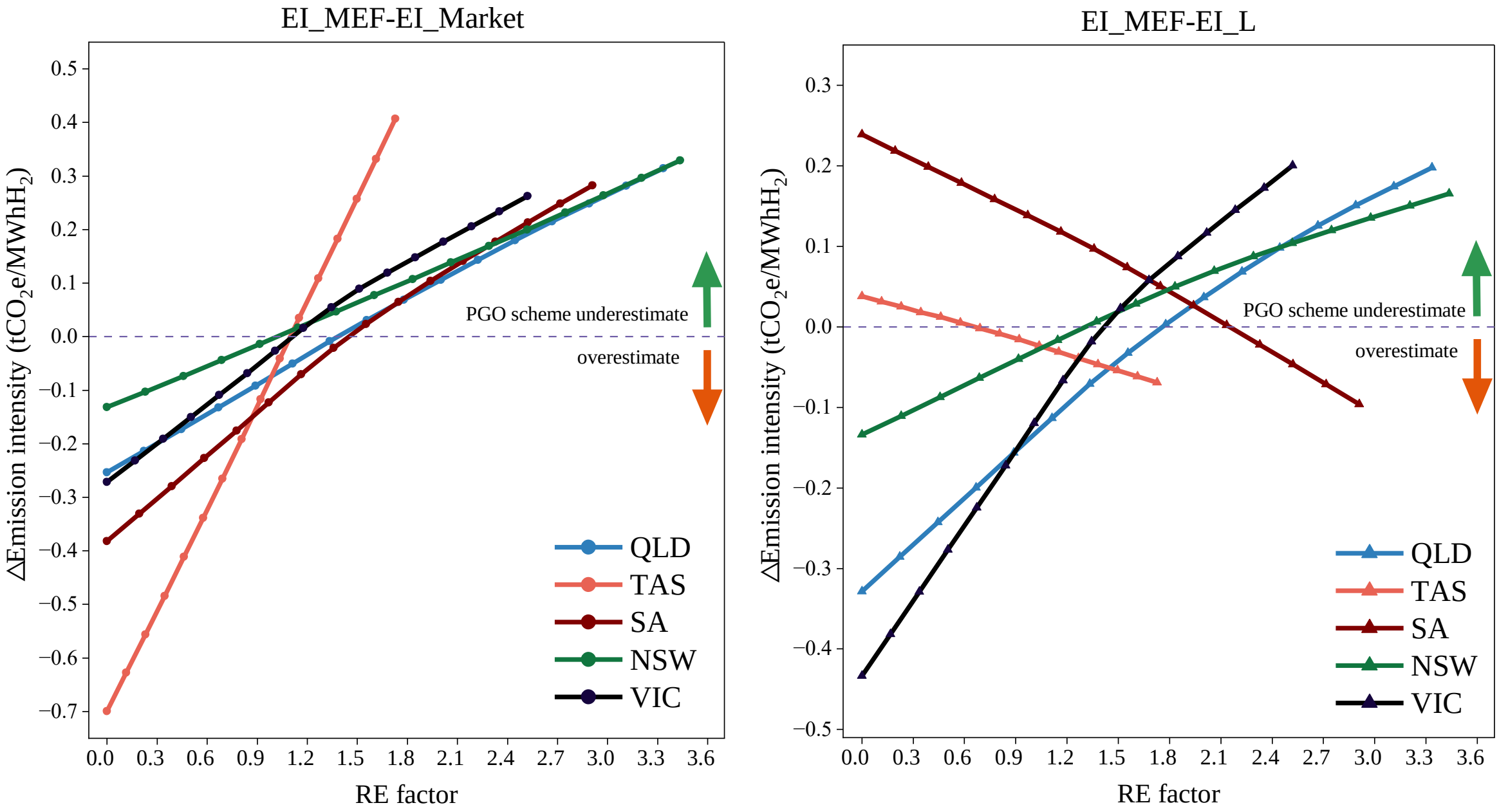

**Figure 16**.The difference in emission intensities calculated by subtracting the PGO market-based method (left) and PGO location-based method (right) from the MEF-based method for different RE factors across the five states in 2023.

## B.2 Verification

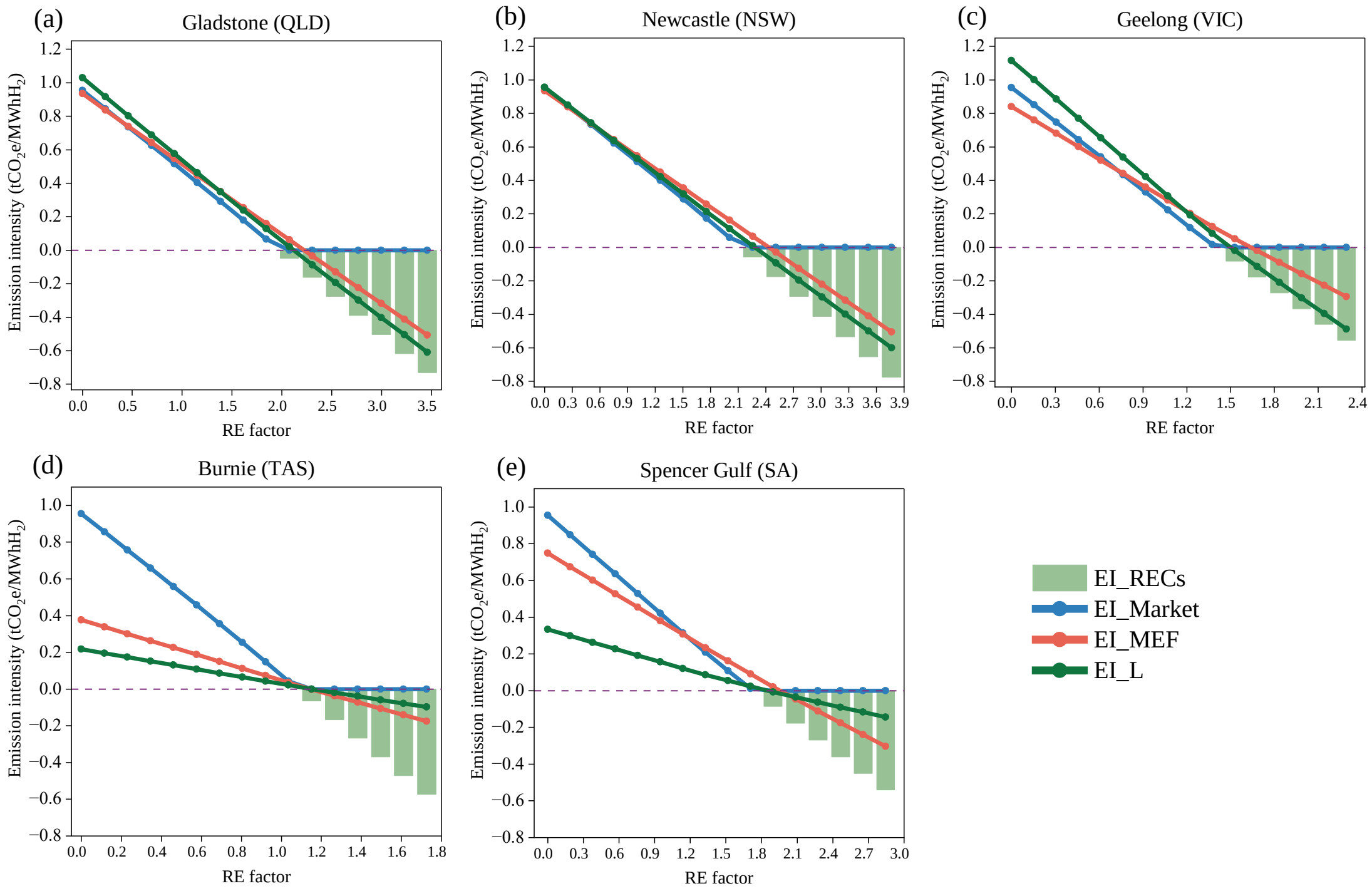

**Figure 17 (a)-(e)** GHG emission intensities of hydrogen under various RE factors in 2021, calculated using three different methods for electricity emissions.

**Figure 17 (a)-(e)** replicate **Figure 2**, utilizing grid data inputs from 2021. These figures illustrate the GHG emission intensities of hydrogen under various RE factors, calculated using three different methods for electricity emissions, along with the differences in emission intensities between the market-based method and the MEF-based method.

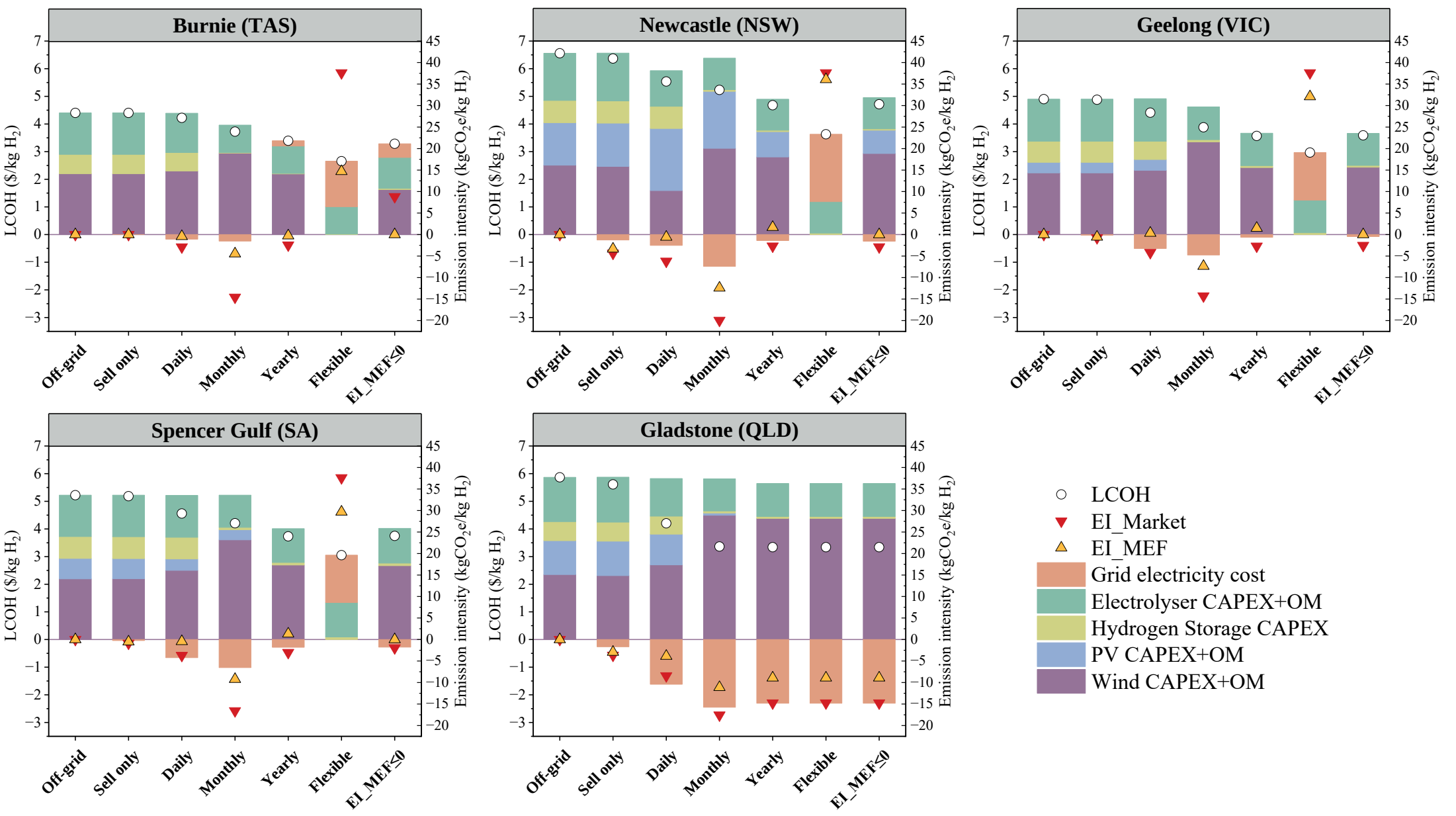

**Figure 18** The LCOH (white dot points) and yearly average emission intensities (red and yellow triangles) of hydrogen calculated for five states in 2021 under different configurations which have been optimized under different scenarios to minimize LCOH.

**Figure 18** depicts the LCOH (white dot points) and yearly average emission intensities (red and yellow triangles) of hydrogen for five states calculated under different configurations optimized to minimize LCOH. **Figure 18** replicates **Figure 3** using the grid data in 2021 for optimization. It is worth noting that the case of Spencer Gulf (SA) in 2021 exhibits different optimization patterns compared to the system optimized using 2023 data, as shown in **Figure 3**. Although the RE capacity factor in Spencer Gulf (SA) is higher in 2021 (0.35) than in 2023 (0.34), indicating more favourable renewable resources, the significantly lower electricity prices in 2021 (59 AUD/MWh) compared to 2023 (92 AUD/MWh) result in the optimized systems under different scenarios in Spencer Gulf (SA) following similar patterns to those observed in TAS, VIC, and NSW. Although the yearly average RE capacity factors of these five cases do not vary significantly (as shown in **Table 6**), the differences in the optimized system results between 2023 (as presented in the **Figure 3**) and 2021 are highly related to variations in electricity prices.

# Conflicts of interest

There are no conflicts to declare.

# Data availability

The data supporting this article have been included in the manuscript.

# Author CRediT statement

Chengzhe Li: Conceptualisation, Investigation, Methodology, Visualization, Writing-Original Draft, Writing-Review and Editing. Lee V. White: Conceptualisation, Investigation, Writing-Review and Editing. Reza Fazeli: Conceptualisation, Investigation, Writing-Review and Editing. Anna Skobeleva: Methodology, Writing-Review and Editing, Michael Thomas: Methodology, Writing-Review and Editing, Shuang Wang: Methodology, Writing-Review and Editing. Fiona J Beck: Conceptualisation, Methodology, Writing-Original Draft, Writing-Review and Editing, Supervision.

# Acknowledgements

This work has been supported by the Heavy Industry Low-carbon Transition Cooperative Research Centre (HILT CRC) whose activities are funded by industry, research and government partners along with the Australian Government's Cooperative Research Centre Program. C.L. acknowledges support from the ANU Centre for Energy Systems.